\documentstyle[12pt]{article}
\begin{document}
\author{S.V. Ludkovsky.}
\title{Non-Archimedean stochastic processes on non-Archimedean manifolds.}
\date{25 October 2002
\thanks{Mathematics subject classification
(1991 Revision) 28C20 and 46S10.} }
\maketitle
\begin{abstract}
Stochastic processes on manifolds over non-Archimedean fields and with
transition measures having values in the field $\bf C$
of complex numbers are defined and investigated.
The analogos of Markov, Poisson and Wiener processes are studied.
For Poisson processes
the non-Archimedean analog of the L\`evy theorem is proved. 
Stochastic antiderivational equations
as well as pseudodifferential equations on manifolds are investigated. 
\end{abstract}
\section{Introduction.} Stochastic processes and stochastic differential
equations on real Banach spaces
and manifolds on them were intensively studied (see, for example,
\cite{beldal,coxmil,dalfo,gihsko,gisk3,ikwat,itkean,malb,mckean,oksb}
and references therein). The considered there stochastic processes
were with values in either real Banach spaces or manifolds on them.
The results of these investigations were used in many mathematical
and theoretical physics problems. In particular stochastic processes
on some Lie groups were studied. On the other hand, the development
of non-Archimedean functional analysis and non-Archimedean quantum
physical theories and quantum mechanics poses problems of developing
measure theory and stochastic processes on non-Archimedean
Banach spaces and manifolds on them \cite{roo,sch1,vla3,khrif,khipr,jang}.
Some steps in this direction were made in \cite{bikvol,evans},
but in these articles only realvalued and complexvalued
stochastic processes were considered.
\par In preceding works of the author measures and stochastic processes
on non-Archimedean Banach spaces and totally disconnected topological
groups with values in non-Archimedean spaces were investigated
\cite{luumnls,lutmf99,lubp2,lulapm,lufpmsp,lustpr}. Quasi-invariant
measures on groups and manifolds were used for investigations
of their representations \cite{luseamb,luseamb2,lubp99}.
\par In this article non-Archimedean stochastic processes and
stochastic antiderivational equations on manifolds on Banach spaces over
non-Archimedean fields are investigated. Moreover,
wider classes of stochastic processes are considered
in this work, than in preceding works of the author.
Analogs of L\`evy processes are studied. More general
classes of analogs of Gaussian measures and Wiener processes
are defined and investigated. Then it is found and proved
a non-Archimedean analog of the It$\hat o$ formula.
\par It is necessary to note that in this article are considered
not only manifolds treated by the rigid geometry, but
much wider classes. For them the existence of an exponential
mapping is proved. A rigid non-Archimedean geometry
serves mainly for needs of the cohomology theory on such manifolds,
but it is too restictive and operates with narrow classes of analytic
functions \cite{freput}. It was introduced at the beginning
of sixties of the 20-th century. Few years later wider classes
of functions were investigated by Schikhof \cite{sch1}.
In this paper classes of functions and antiderivation operators
by Schikhof and their generalizations from works \cite{lutmf99,lubp2}
are used.
\par The results of this paper permit to consider stochastic
processes on non-Archimedean manifolds as well as more
general classes of stochastic processes on non-Archimedean
Banach spaces and totally disconnected topological groups.
Some other principal differences of the classical and
non-Archimedean stochastic analyses are discussed in \cite{lustpr}.
\section{Stochastic processes for non-Archimedean
locally $\bf K$-convex spaces.}
To avoid misunderstanding we first give our definitions and notations.
\par {\bf 2.1. Definitions and Notes.}
A measurable space $(\Omega ,{\sf F})$ with a probability realvalued
$\sigma $-additive measure $\lambda $ on a covering $\sigma $-algebra
$\sf F$ of a set
$\Omega $ is called a probability space and it is denoted by
$(\Omega ,{\sf F},\lambda )$.
In the case of a complexvalued $\sigma $-additive
measure $\lambda $ we suppose, that its variation $|\lambda |$
is a probability realvalued $\sigma $-additive measure,
which is natural, since $|\lambda |$ is a nonegative
$\sigma $-additive measure.
Points $\omega \in \Omega $ are called
elementary events and values $\lambda (S)$  
probabilities of events $S\in \sf F$. A measurable map
$\xi : (\Omega ,{\sf F})\to (X,{\sf B})$ is called a random variable
with values in $X$, where ${\sf B}$ is a covering $\sigma $-algebra
of subsets of a locally $\bf K$-convex space $X$,
$\xi ^{-1}({\sf B})\subset \sf F$,
where $\bf K$ is a non-Archimedean field complete as an ultrametric space. 
\par The random variable $\xi $ induces a normalized measure
$\nu _{\xi }(A):=
\lambda (\xi ^{-1}(A))$ in $X$ and a new probability space 
$(X,{\sf B},\nu _{\xi }).$
\par Let $T$ be a set with a covering $\sigma $-algebra $\cal R$ and a
$\sigma $-additive measure
$\eta : {\cal R}\to \bf C$. Consider the following Banach space
$L^q(T,{\cal R},\eta ,H)$ as the completion of the set
of all ${\cal R}$-step functions $f: T\to H$ relative to the following 
norm:
\par $(1)\quad \| f\|_{\eta ,q}:=[\int_T \| f(t)\|_H^q
|\eta |(dt)]^{1/q}$ for $1\le q<\infty $ and
\par $(2)\quad \| f\|_{\eta ,\infty }:=ess-\sup_{t\in T,\eta }
\| f(t)\|_H$, where $H$ is a Banach space over $\bf K$,
$|\eta |$ is a variation of $\eta $, that is, $|\eta |$
is a realvalued $\sigma $-additive measure.
For $0<q<1$ this is the metric space with the metric
\par $(3)\quad \rho _q(f,g):=[\int_T\| f(t)-g(t)\|_H^q |\eta |(dt)]^{1/q}.$  
\par Consider now a complete locally $\bf K$-convex space $H$,
then $H$ is a projective limit of Banach spaces
$H=\lim \{ H_{\alpha },\pi ^{\alpha }_{\beta }, \Upsilon \} $,
where $\Upsilon  $ is a directed set, $\pi ^{\alpha }_{\beta }:
H_{\alpha }\to H_{\beta }$ is a $\bf K$-linear continuous mapping
for each $\alpha \ge \beta $, $\pi _{\alpha }: H\to H_{\alpha }$
is a $\bf K$-linear continuous mapping such that $\pi ^{\alpha }_{\beta }
\circ \pi _{\alpha }=\pi _{\beta }$ for each $\alpha \ge \beta $
(see \S 6.205 \cite{nari}).
Each norm $p_{\alpha }$ on $H_{\alpha }$ induces a seminorm
${\tilde p}_{\alpha }$ on $H$. If $f: T\to H$, then $\pi _{\alpha }\circ
f=:f_{\alpha }: T\to H_{\alpha }$. In this case $L^q(T,{\cal R},\eta ,H)$
is defined as a completion of a family of all step functions
$f: T\to H$ relative to the family of seminorms
\par $(1')\quad \| f\|_{\eta ,q,\alpha }:=[\int_T{\tilde p}_{\alpha }
(f(t))^q|\eta |(dt)]^{1/q}$, $\alpha \in \Upsilon $, for $1\le q<\infty $ and
\par $(2')\quad \| f\|_{\eta ,\infty ,\alpha  }:=ess-\sup_{t\in T,\eta }
{\tilde p}_{\alpha }(f(t))$, $\alpha \in \Upsilon $,
or pseudometrics
\par $(3')\quad \rho _{q,\alpha }(f,g):=[\int_T{\tilde p}_{\alpha }
((f(t)-g(t))^q|\eta |(dt)]^{1/q}$, $\alpha \in \Upsilon $, for $0<q<1$.
Consequently, $L^q(T,{\cal R},\eta ,H)=$
$\lim \{ L^q(T,{\cal R},\eta ,H_{\alpha }),\pi ^{\alpha }_{\beta },
\Upsilon \} $. For example, $T$ may be a subset of $\bf F$ or of $\bf R$,
where $\bf F$ denotes a non-Archimedean field.
\par If $T$ is a zero-dimensional $T_1$-space, then denote by
$C^0_b(T,H)$ the Banach space of all continuous bounded functions
$f: T\to H$ supplied with the norm:
\par $(4)\quad \| f\|_{C^0}:=\sup_{t\in T} \| f(t)\|_H<\infty $. \\
\par For a set $T$ and a complete locally $\bf K$-convex 
space $H$ over $\bf K$ consider the 
product $H^T$ of $\bf K$-convex spaces
$H^T:=\prod_{t\in T}H_t$ in the product topology,
where $H_t:=H$ for each $t\in T$.
\par Then take on either $X:=X(T,H)=L^q(T,{\cal R},\eta ,H)$ or $X:=
X(T,H)=C^0_b(T,H)$ or on $X=X(T,H)=H^T$ a covering $\sigma $-algebra
${\sf B}$, for example, ${\sf B}\supset Bf (X)$,
where $Bf (X)$ denotes the Borel $\sigma $-algebra of subsets of $X$
for a given topology on $X$.
Consider a random variable
$\xi : \omega \mapsto \xi (t,\omega )$ with values in $(X,{\sf B})$,
where $t\in T$.
\par Events $S_1,...,S_n$ are called independent in total if
$P(\prod_{k=1}^nS_k)=\prod_{k=1}^nP(S_k)$. Subalgebras
${\sf F}_k\subset {\sf F}$ are said to be independent if
all collections of events $S_k\in {\sf F}_k$ are independent in total, 
where $k=1,...,n$, $n\in \bf N$. To each collection of random variables
$\xi _{\gamma }$ on $(\Omega ,{\sf F})$ with $\gamma \in \Upsilon $
is related the minimal $\sigma $-algebra ${\sf F}_{\Upsilon }\subset \sf F$
with respect to which all $\xi _{\gamma }$ are measurable, where $\Upsilon $
is a set.
Collections $\{ \xi _{\gamma }: $ $\gamma \in \Upsilon _j \} $
are called independent if such are 
${\sf F}_{\Upsilon _j}$, where $\Upsilon _j\subset \Upsilon $ for each
$j=1,...,n,$ $n\in \bf N$.
\par Let $T$ be such that $card(T)>n$. For $X=C^0_b(T,H)$ or $X=H^T$
define $X(T,H;(t_1,...,t_n);(z_1,...,z_n))$ as a closed submanifold
in $X$ of all $f: T\to H$, $f\in X$ such that $f(t_1)=z_1,...,f(t_n)=z_n$,
where $t_1,...,t_n$ are pairwise distinct points in $T$ and
$z_1,...,z_n$ are points in $H$. For $X=L^q(T,{\cal R},\eta ,H)$ and 
pairwise distinct points $t_1,...,t_n$ in $T\cap supp (|\eta |)$
define a set
$X(T,H;(t_1,...,t_n);(z_1,...,z_n))$ as a closed submanifold
which is the completion relative to the
metric $\rho _q$ (or a family of pseudometrics $ \{ \rho _{q,\alpha }:
\alpha \} $ respectively), where $0<q\le \infty $,
of a family of $\cal R$-step functions $f: T\to H$ such that
$f(t_1)=z_1,...,f(t_n)=z_n$. In these cases 
$X(T,H;(t_1,...,t_n);(0,...,0))$ is the proper $\bf K$-linear subspace 
of $X(T,H)$ such that $X(T,H)$ is isomorphic with
$X(T,H;(t_1,...,t_n);(0,...,0))\oplus H^n$, since if
$f\in X$, then $f(t)-f(t_1)=:g(t)\in X(T,H;t_1;0)$
(in the third case we use that $T\in \cal R$ and hence there exists
the embedding $H\hookrightarrow X$). For $n=1$ and $t_0\in T$
and $z_1=0$ we denote $X_0:=X_0(T,H):=X(T,H;t_0;0)$.
\par {\bf 2.2. Definition.} We define a stochastic process
$\xi (t,\omega )$ with values in $H$ as a random variable such that:
\par $(i)$ the differences $\xi (t_4,\omega )-\xi (t_3,\omega )$ 
and $\xi (t_2,\omega )-\xi (t_1,\omega )$ are independent
for each chosen $(t_1,t_2)$ and $(t_3,t_4)$ with $t_1\ne t_2$,
$t_3\ne t_4$, such that either $t_1$ or $t_2$ is not in the two-element set
$ \{ t_3,t_4 \} ,$ when $T$ is a subset in $\bf R$ we suppose additionally,
that $t_1<t_2$ and $t_3<t_4$, where $\omega \in \Omega ;$
\par $(ii)$ the random variable $\xi (t,\omega )-\xi (u,\omega )$ has 
a distribution $\mu ^{F_{t,u}},$ where $\mu $ is a probability 
complexvalued measure on $(X(T,H),{\sf B})$ from \S 2.1, 
$\mu ^g(A):=\mu (g^{-1}(A))$ for $g: X\to H$ such that
$g^{-1}({\cal R}_H)\subset \sf B$ and for each $A\in
{\cal R}_H$, a continuous linear operator $F_{t,u}: X\to H$ is given by
the formula $F_{t,u}(\xi ):=\xi (t,\omega )-\xi (u,\omega )$ 
for each $\xi \in L^r(\Omega ,{\sf F},\lambda ;X),$
where $0< r\le \infty ,$
${\cal R}_H$ is a covering $\sigma $-algebra of $H$ such that
$F_{t,u}^{-1}({\cal R}_H)\subset \sf B$ for each $t\ne u$ in $T$;
\par $(iii)$ we also put $\xi (0,\omega )=0,$  
that is, we consider a $\bf K$-linear subspace
$L^r(\Omega ,{\sf F},\lambda ;X_0)$ of $L^r(\Omega ,{\sf F},\lambda ;X)$,
where $\Omega \ne \emptyset $, $X_0$
is the closed subspace of $X$ as in \S 2.1.
\par It is seen that $\xi (t,\omega )$ is a Markov process with 
a transition measure $P(u,x,t,A)=\mu ^{F_{t,u}}(A-x)$.
\par This definiton is justified by the following Theorem.
\par {\bf 2.3.} {\bf Theorem.} {\it 
Let either $X=C^0_b(T,H)$ or $X=H^T$ or $X=L^g(T,{\cal R},\eta ,H)$ with
$0<g \le \infty $ be the same spaces as in \S 2.1,
where the valuation group $\Gamma _{\bf K}$ is discrete
in $(0,\infty )$.
Then there exists a family $\Psi $
of pairwise inequivalent stochastic processes on 
$X$ of a cardinality $card (\Psi )\ge card (T) card (H)$
or $card (\Psi )\ge card ({\cal R}) card (H)$ respectively.}
\par {\bf The proof} is analogous to that of
Theorem I.4.3 \cite{lustpr} and Theorem 4.3 \cite{lukhr}
as well as Note 4.4 \cite{lukhr} can be applied to the considered
here case also.
\par {\bf 2.4. Definition.} Let $T$ be an additive group contained
in $\bf R$.
Consider a stochastic process $\xi \in L^r(\Omega ,{\sf F},\lambda ,
X_0(T,H))$ such that a transition measure
has the form 
\par $P(t_1,x,t_2,A):=P(t_2-t_1,x,A):=\exp (-\rho (t_2-t_1))
P(A-x)$ \\
(see \S 3.2 \cite{lukhr} and \S 2.2 above)
for each $x\in H$ and $A\in {\cal R}_H$
and $t_1$ and $t_2$ in $T$, where $\rho >0$ is a constant. 
Then such process is called the Poisson process.
\par {\bf 2.5. Proposition.} {\it Let
\par $P(A-x)=\int_HP(-x+dy)P(A-y)$ 
for each $x\in H$ and $A\in {\cal R}_H$, where $P$ is a probability
measure on $H$ and $T$ is an interval in $\bf R$,
\par then there exists
a measure $\mu $ on $X_0(T,H)$ for which the Poisson process
exists.}
\par {\bf Proof.} We take 
\par $\mu _{t_1,...,t_n}:=P(t_2-t_1,0,*)...P(t_n-t_{n-1},0,*)$ on 
\par ${\cal R}_{t_1,...,t_n}:={\cal R}_{t_1}\times ...\times {\cal R}_{t_n}$\\
for each pairwise distinct points $t_1,...,t_n\in T$, where 
\par $\mu _{t_1,...,t_n}=\pi _{t_1,...,t_n}(\mu )$ and 
\par $\pi _{t_1,...,t_n}: X_0(T,H)\to H_{t_1}\times ...\times H_{t_n}$ \\
is the natural projection, $H_t=H$ for each $t\in T$.
There is a family $\Lambda $ of all finite subsets of $T$
directed by inclusion.
In view of the Kolmogorov theorem (see Theorems I.1.3, I.1.4
\cite{dalfo} and \cite{kolfp}) the cylindrical distribution
$\mu $ generated by the family $P(t_2-t_1,x,A)$ has an extension
to a measure on $X_0(T,H)$. All others conditions are satisfied 
in accordance with \S 2.2 and \S 2.4.
\par {\bf 2.6. Remark.}
Let $K$ be a complete uniform (or in particular ultrauniform)
Tychonoff space. Put 
\par ${\tilde K}^n:=(x\in K^n: x_i\ne x_j$
for each $i\ne j)$. \\
Supply ${\tilde K}^n$ with a product topology.
Let also $B^n_K$ denotes the collection of all
$n$-point subsets of $K$.
For each subset $A\subset K$ a number mapping $N_A: B^n_K\to \bf N_o$
is defined by the following formula: $N_A(\gamma ):=card(\gamma \cap A)$,
where ${\bf N}:=\{ 1,2,3,... \} $, ${\bf N_o}:=\{ 0,1,2,3,... \} $.
\par {\bf 2.7. Definitions and Remarks.} As usually let 
\par $B_K:=\bigoplus_{n=0}^{\infty }B^n_K$, \\
where $B^0_K:=\{ \emptyset \} $ is a singleton,
$B_K\ni x=(x_n: x_n\in B^n_K, n=0,1,2,... )$. 
If a complete (ultra)uniform space $X$ is not compact,
then there exists an increasing sequence of subsets $K_n\subset X$
such that $X=\bigcup_nK_n$ and $K_n$ are complete spaces in the
uniformity induced from $X$. Moreover, $K_n$ can be chosen clopen in $X$,
when $X$ is ultrauniform.
Then the following space 
\par $\Gamma _X:=
\{ \gamma : \gamma \subset X$ and $card(\gamma \cap K_n)<\infty $
for each $n \} $ \\
is called the configuration space and it is isomorphic
with the projective limit $pr-\lim \{ B_{K_n}, \pi ^n_m, {\bf N} \} $,
where $\pi ^n_m(\gamma _m)=\gamma _n$ for each $m>n$ and $\gamma _n
\in B_{K_n}$.
Then $\prod_{n=1}^{\infty }B_{K_n}=:Y$ in the Tychonoff product
topology is ultrauniformizable, that induces the ultrauniformity in
$\Gamma _X$, for example, 
\par $\rho (x,y):=d_n(x_n,y_n)p^{-n}$ is the pseudoultrametric
in $\Gamma _X$ for a family of pseudoultrametrics $d_n$ on $B_{K_n}$, \\
where $n=n(x,y):=\min_{(x_j\ne y_j)}j$,
$x=(x_j: j\in {\bf N}, x_j\in B_{K_j} ) $, $1<p\in \bf N$,
since each complete ultrauniform space is a projective limit
of complete ultrametric spaces \cite{eng,lufpm}
(about $d_n$ in the case of ultrametric spaces see \cite{luseamb2}).
\par  Let $K\in \{ K_n: n\in {\bf N} \} $, then $m_K$
denotes the restriction $m|_K$, where $m: {\cal R}\to \bf C$
is a $\sigma $-additive measure on a
covering $\sigma $-algebra ${\cal R}_m$ of $X$, $K_n\in {\cal R}_m$
for each $n\in \bf N$. Suppose that $K^n_l\in {\cal R}_{m^n}$
for each $n$ and $l$ in $\bf N$, where ${\cal R}_{m^n}$ is the completion
of the covering $\sigma $-algebra ${\cal R}^n$ of $X^n$
relative to the product measure
$m^n=\bigotimes_{j=1}^nm_j$, $m_j=m$ for each $j$.
Then $m^n_K:=\bigotimes_{j=1}^n(m_K)_j$ is a measure
on $K^n$ and hence on ${\tilde K}^n$, when $m$ is 
such that $|m|(K^n\setminus {\tilde K}^n)=0$,
where $(m_K)_j=m_K$ for each $j$.
Let $|m|(X) <\infty $. Then 
\par $(i)$ $P_{K,m}:=\exp (-m(K))\sum_{n=0}
^{\infty }m_{K,n}/n!$ \\
is a measure on ${\cal R}(B_K)$, where
\par ${\cal R}(B_K)=B_K\cap (\bigoplus_{n=0}^{\infty }{\cal R}_{m^n}),$ \\
$m_{K,0}$ is a probability measure on the singleton $B^0_K$, and
$m_{K,n}$ are images of $m^n_K$ under the following mappings:
\par $p^n_K: {\tilde K}^n\ni (x_1,...,x_n)\mapsto
\{ x_1,...,x_n \} \in B^n_K$. \\
Such system of measures $P_{K,n}$ is consistent, that is, 
\par $\pi ^n_l(P_{K_l,m})=P_{K_n,m}$ for each
$n\le l$. \\
This defines the unique measure $P_m$ on 
${\cal R}(\Gamma _X)$, which is called the Poisson measure, where
$\pi _n: Y\to B_{K_n}$ is the natural projection for each
$n\in \bf N$.
For each $n_1,..,n_l\in \bf N_o$ and disjoint subsets
$B_1,...,B_l$ in $X$ belonging to ${\cal R}_m$
there is the following equality:
\par $(ii)$ $P_m(\bigcap_{j=1}^l \{ \gamma : card(\gamma \cap B_i)=n_i \}
)=\prod_{i=1}^lm(B_i)^{n_i} \exp(-m(B_i))/n_i!$.
\par There exists the following embedding $\Gamma _X\hookrightarrow
S_X$, where 
\par $S_X:=\lim \{ E_{K_n}, \pi ^n_m, {\bf N} \} $ 
is the limit of an inverse mapping sequence,
\par $E_K:=\bigoplus_{l=0}
^{\infty }K^l$ for each $K\in \{ K_n: n=0,1,2,... \} $. \\
The Poisson measure $P_m$ on ${\cal R}(\Gamma _X)$ considered
above has an extension on ${\cal R}(S_X)$ such that
$|P_m|({S_X\setminus \Gamma _X})=0$. If each $K_n$ is a 
complete $\bf K$-linear space (not open in $K$),
then $E_K$ and $S_X$ are complete $\bf K$-linear spaces, since 
\par $S_X\subset (\prod_{n=1}^{\infty }
E_{K_n})$. \\
Then on ${\cal R}(S_X)$ there exists a Poisson measure
$P_m$, but without the restriction 
$|m^n_K|({K^n\setminus {\tilde K}^n})=0$, where
\par $(iii)$ $P_{K,m}:=\exp (-m(K))\sum_{n=0}^{\infty }
m^n_K/n!$, 
\par $\pi ^n_l(P_{K_l,m})=P_{K_n,m}$ for each $n\le l$.
\par {\bf 2.8. Corollary.} {\it Let suppositions of Proposition
2.5 be satisfied with $H=S_X$ for a complete $\bf K$-linear space
$X$ and $P(A)=P_m(A)$ for each $A\in {\cal R}(S_X)$
(see \S 2.7), then there exists
a measure $\mu $ on $X_0(T,H)$ for which the Poisson process
exists.}
\par {\bf 2.9. Definition.} The stochastic process of Corollary 2.8
is called the Poisson process with values in $X$.
\par {\bf 2.10. Note.} If $\xi \in L^r(\Omega , {\sf F}, \lambda ;
X_0(T,Y))$ is a stochastic process with values in a Hilbert space
$Y$ over $\bf C$, then its mean value
for $t_1, t_2\in T$ is defined by the following formula:
$$(i)\quad M_{t_1,t_2}(\eta ):= \int_Yy P(t_1,0,t_2,dy),$$
where $P(t_1,y_1,t_2,A)$ is a transition probability of $\xi $
corresponding to $\xi (t_1,\omega )=y_1$, $\xi (t_2,\omega )\in A$,
$A\in Bf (Y)$. For $t_1=t_0$ we may simply write
$M_{t_2}$, if $t_1$ and $t_2$ are definite moments, then
they may be omitted and we may write $M$ instead of $M_{t_1,t_2}$.
Let $H=\bf K$ be a field, where ${\bf Q_p}\subset {\bf K}\subset \bf C_p$,
let also $\lambda $ be a probability realvalued measure.
Let $T$ be an interval $[t_0,R]$ in $\bf R$, where $R>t_0$. 
Consider a multiplicative character for a field $\bf K$,
$\pi : {\bf K}\setminus \{ 0 \} \to \bf C$,
$\pi =\pi _a$ for some $a\in \bf C$ such that 
$\pi _a(x):=|x|_{\bf K}^{a-1}\pi _0(x|x|_{\bf K})$,
where $\pi _0: S_1\to S^1$ is a multiplicative character
on $S_1:= \{ x\in {\bf K}: |x|_{\bf K}=1 \} $,
${\bf C_p}\supset {\bf K}\supset \bf Q_p$,
$S^1:=\{ z\in {\bf C}: |z|=1 \} $ (about a character see, for example,
\S VI.25 \cite{hew} and \S III.2 \cite{vla3}).
We take $a$ with $Re (a)>1$ and consider an extension of
$\pi _a$ as a continuous function such that $\pi _a(0)=0$.
\par {\bf 2.11. Theorem.} {\it Let $\psi $ be a
continuously differentiable function,
from an interval $T\subset \bf R$ into $\bf R$ and $\psi (0)=0$.
Then there exists a stochastic process $\xi (t,\omega )$ such that 
\par $M_t(\exp(-\rho \pi [\xi (t,\omega )]))=\exp (-t\psi (\rho ))$ \\
for each $t$ in $T$ and each constant $\rho >0$, where
$\pi : {\bf K}\setminus \{ 0 \} \to {\bf C}$
is a multiplicative continuous character as in \S 2.10.}
\par {\bf Proof.} We consider solution of the following equation 
$$M_t[\exp (-\rho \pi [\xi (t,\omega )])]=\exp (-t\psi (\rho ))$$
taking $t_0=0$ without loss of generality.
Then $e(t)=e(t-s)e(s)$ 
for each $t$ and $s$ in $T$ and each $\rho >0$, where 
\par $e_{\rho }(t):=e(t):=M_t(\exp (-\rho \pi [\xi (t,\omega )])).$ \\
Hence 
\par $\partial e_{\rho }(t)/\partial \rho =-t\psi '(\rho )\exp (-t\psi
(\rho ))$, consequently, 
\par $\psi '(\rho )=t^{-1}(\int_{\bf K}\pi (l)
\exp (-\rho \pi (l)) P(\{ \omega : \xi (t,\omega )\in dl \} )$ \\
for each $t\ne 0$. In particular,
\par $\psi '(\rho )=\lim_{t\to 0, t\ne 0}
t^{-1}(\int_{\bf K}\pi (l) \exp (-\rho \pi (l))
P(\{ \omega : \xi (t,\omega ) \in dl \} )$. \\
By the conditions of this theorem we have
\par  $\psi (\alpha )=\int_0^{\alpha } \psi '(\beta )d\beta $. \\
Consider a $\sigma $-additive measure $m$ on a separating covering ring
${\cal R}({\bf K})$ such that ${\cal R}({\bf K})\supset
Bf ({\bf K})\cup \{ 0 \} $ with values in $\bf C$ given by the formula
\par $m(dl):=\lim_{t\to 0, t\ne 0} \pi (l) P(\{ \omega : \xi (t,\omega )
\in dl \} )/t$ \\
on ${\bf K}\setminus \{ 0 \} $, $m ( \{ 0 \} )=m_0$ and consider
a $\sigma $-additive measure $n$ such that 
\par $m(dl)=\pi (l) n(dl)$ for $l\ne 0$. \\
Therefore, 
\par $\psi (\rho )=\int_0^{\rho }(\int_{\bf K}\exp
(-\beta \pi (l))m(dl)) d\beta .$ \\
The character $\pi $ is continuous and multiplicative, that is,
$\pi (ab)=\pi (a)\pi (b)$ for each $a$ and $b\in {\bf K}\setminus \{ 0 \} $.
From $\psi (0)=0$ we have $e_0(t)=1$ for each $t\ge 0$, consequently,
\par $\psi (\rho )=\rho m_0+\int_{{\bf K}\setminus \{ 0 \}}
[1-\exp (-\rho \pi (l))] n(dl)$, \\
since $\pi $ as the continuous function has the extension on $\bf K$
with $\pi (0)=0$ and
\par  $\lim_{l \to 0, l \ne 0}
[1-\exp (-\rho \pi (l))]/\pi (l)=\rho $ and 
$$\lim_{\rho \to 0, \rho \ne 0}
\int_{B({\bf K},0,k)}[1-\exp (-\rho \pi (l))] n(dl)=0$$
for each $k>0$. Then
\par $\psi (\rho )=\int_{\bf K}[1-\exp (-\rho \pi (l))] n(dl)$.
We search a solution of the problem in the form
\par $\pi (\xi (t,\omega ))=tm_0+
\int_{\bf K}\pi (l) {\cal \eta }([0,t],dl,\omega )$, \\
where ${\cal \eta }(dt, dl, \omega )$ is the realvalued
$\sigma $-additive measure on $Bf (T)\times {\cal R}({\bf K})$
for each $\omega \in \Omega $ such that
its moments satisfy the Poisson distribution with the Poisson measure
$P_{tn}$, that is, 
\par $M_t[{\cal \eta }^k([0,t],dl,\omega )]=
\sum_{s\le k}a_{s,k}(tn)^s(dl)/s!$ \\
for each $0<t\in T$, where
$$a_{k,j}=\sum_{s_1+...+s_k=j, s_1\ge 1,...,s_k\ge 1} j!/(s_1!...s_k!)$$
for each $k\le j.$
Using the fact that the set of step functions is dense
in $L^r({\bf K},{\cal R}({\bf K}),n,{\bf C_p})$
we get
$$M_t[\exp (-\rho \int_{\bf K}\pi (l){\cal \eta }([0,t],dl,\omega ))]
=\lim_{\cal Z}M_t[\prod_j\exp (-\rho \pi (l_j)
{\cal \eta }([0,t],\delta _j,\omega ))]$$
$$=\lim_{\cal Z}\prod_jM_t[\exp (-\rho \pi(l_j)
{\cal \eta }([0,t],\delta _j,\omega ))]
=\lim_{\cal Z}\exp (-\rho t\sum_j(1-\exp (-\rho \pi (l_j)))n(\delta _j))$$
$$=\exp [-\rho t\int_{\bf K} \{ 1-\exp (-\rho \pi (l)) \} n(dl)],$$
where ${\cal Z}$ is an ordered
family of partitions ${\cal U}$ of $\bf K$ into disjoint union of elements
of ${\cal R}({\bf  K})$, ${\cal U}\le \cal V$ in ${\cal Z}$
if and only if each element of the disjoint covering ${\cal U}$
is a union of elements of $\cal V$, $l_j\in \delta _j\in {\cal U}\in \cal Z$.
We get the equation 
$$M_t[\exp (-\rho \pi [\xi (t,\omega )])]=
M_t[\exp (-\rho \int_{\bf K}\pi (l){\cal \eta }([0,t],dl,\omega ))].$$
This defines the stochastic process $\pi [\xi (t,\omega )]$
with the probability space $(\Omega ,{\sf F},\lambda )$.
If $f\in L^r(\Omega ,{\sf F},\lambda ;X_0(T,{\bf C}))$
and $f(\Omega \times T)\subset \pi ({\bf K})\subset \bf C$,
then there exists $h\in L^r(\Omega ,{\sf F},\lambda ;
X_0(T,{\bf K}))$ such that $f=\pi [h]$.
Since $\pi $ is continuous on $\bf K$ and locally constant on ${\bf K}
\setminus \{ 0 \} $, then
due to \S 2.3 and Corollary 2.8 above there exists a
$\bf K$-valued stochastic process
$\xi $ for a given complexvalued stochastic process
$\pi [\xi ]$ with the measure space $(\Omega ,{\sf F},\lambda )$
(see also \S 4.3 \cite{lukhr}).
\par {\bf 2.12.1. Notes.} For a continuous function
$c(t)$ and a nonnegative measure $n$ on $T\times \bf K$
such that $\int_{{\bf K}\setminus \{ 0 \} } (1-\exp
(-\pi (l))n((0,t]\times dl)<\infty $ for each $t>0$
there exists a $\bf K$-valued stochastic process
$\xi (t,\omega )$ such that
$$M_{t_1,t_2}[\exp \{ -\rho (\pi (\xi (t_2,\omega ))-
\pi (\xi (t_1,\omega ))) \} ]=\exp [- \rho (c(t_2)-c(t_1))]$$
$$-\int_{{\bf K}\setminus \{ 0 \} }[1-\exp (-\rho \pi (l))]
n((t_1,t_2]\times dl)]$$ for each $\rho >0$ and
$t_1< t_2$ with a Poisson measure $\eta $ having
mean $n(dt\times dl)$. This can be proved analogously
to \S 4.10 \cite{itkean} using \S 2.11 above and with the help
of disjoint pavings of $\bf K$ by clopen balls instead of
intervals $(l_1,l_2]$ in the real case. For this put
$$\pi (\xi (t,\omega ))=c(t)+\int_{{\bf K}\setminus \{ 0 \} }
\pi (l)\eta ([0,t]\times dl)$$
for each $t\ge 0$, where $\eta (dt\times dl)$ is a number
of jumps of magnitude $s\in dl$ in time $dt$.
\par {\bf 2.12.2. Notes.}
In \cite{lubp2} and in II \cite{lustpr} was considered
an analog of a Gaussian measure and of a Wiener process.
That construction is generalized below and additional properties
are proved conserning moments of a Gaussian measure and
an analog of the It$\hat o$ formula.
\par Let $X$ be a locally $\bf K$-convex space equal to
a projective limit $\lim \{ X_j, \phi ^j_l, \Upsilon \} $
of Banach spaces over a local field $\bf K$ such that
$X_j=c_0(\alpha _j,{\bf K})$, where the latter space consists
of vectors $x=(x_k: k\in \alpha _j )$, $x_k\in \bf K$,
$\| x \| :=\sup_k |x_k|_{\bf K}<\infty $ and such that
for each $\epsilon >0$ the set $\{ k: |x_k|_{\bf K}>\epsilon \} $
is finite, $\alpha _j$ is a set, that is convenient
to consider as an ordinal due to Kuratowski-Zorn lemma \cite{eng,roo};
$\Upsilon $ is an ordered set,
$\phi ^j_l: X_j\to X_l$ is a $\bf K$-linear continuous mapping
for each $j\ge l\in \Upsilon $, $\phi _j: X\to X_j$ is a projection
on $X_j$, $\phi _l\circ \phi ^j_l=\phi _j$ for each $j\ge l\in \Upsilon $,
$\phi ^l_k\circ \phi ^j_l=\phi ^j_k$ for each $j\ge l\ge k$ in $\Upsilon $. 
Consider also a locally $\bf R$-convex space, that is
a projective limit $Y=\lim \{ l_2(\alpha _j,{\bf R}),
\psi ^j_l, \Upsilon \} $, where $l_2(\alpha _j,{\bf R})$ is the real Hilbert
space of the topological weight $w(l_2(\alpha _j,{\bf R}))=
card (\alpha _j)\aleph _0$. Suppose $B$ is a symmetric
nonegative definite (bilinear) nonzero functional $B: Y^2\to \bf R$.
\par {\bf 2.13. Definitions and Notes.}
A measure $\mu =\mu _{q,B,\gamma }$ on $X$
with values in $\bf R$ is called a $q$-Gaussian measure, if
its characteristic functional $\hat \mu $ has the form
$$\hat \mu (z)=\exp [-B(v_q(z),v_q(z))]\chi _{\gamma }(z)$$
on a dense $\bf K$-linear subspace ${\sf D}_{q,B,X}$ in $X^*$ of all
continuous $\bf  K$-linear functionals $z: X\to \bf K$
of the form $z(x)=z_j(\phi _j(x))$ for each $x\in X$
with $v_q(z)\in {\sf D}_{B,Y}$,
where $B$ is a nonegative definite symmetric operator
(that is, bilinear $\bf R$-valued symmetric functional)
on a dense $\bf R$-linear subspace ${\sf D}_{B,Y}$ in $Y^*$,
$B: {\sf D}_{B,Y}^2\to \bf R$,
$j\in \Upsilon $ may depend on $z$, $z_j: X_j\to \bf K$ is a continuous
$\bf K$-linear functional such that $z_j=\sum_{k\in \alpha _j}
e^k_jz_{k,j}$ is a countable convergent series such that
$z_{k,j}\in \bf K$, $e^k_j$ is a continuous $\bf K$-linear functional
on $X_j$ such that $e^k_j(e_{l,j})=\delta ^k_l$ is the Kroneker delta
symbol, $e_{l,j}$ is the standard orthonormal (in the non-Archimedean sence)
basis in $c_0(\alpha _j,{\bf K})$, $v_q(z)=v_q(z_j):=\{
|z_{k,j}|_{\bf K}^{q/2}: k\in \alpha _j \} $.
It is supposed that $z$ is such that $v_q(z)\in l_2(\alpha _j,{\bf R})$,
where $q$ is a positive constant, $\chi _{\gamma }(z): X\to S^1$
is a continuous character such that $\chi _{\gamma }(z)=
\chi (z(\gamma ))$, $\gamma \in X$, $\chi : {\bf K}\to S^1$
is a character of $\bf K$ as an additive group
(about a character see,
for example, \S VI.25 \cite{hew} and \S III.1 \cite{vla3}).
\par A symmetric nonegative definite operator $B$ is called
a correlation operator of a measure $\mu $. If $Y$ is a Hilbert space
with a scalar product $(*,*)$, then due to the Riesz theorem
there exists $E\in L(Y)$ such that $B(y_1,y_2)=(Ey_1,y_2)$
for each $y_1,$ $y_2\in Y$.
Therefore, $B$ is also called operator.
\par {\bf 2.14. Proposition.} {\it A $q$-Gaussian measure
on $X$ is $\sigma $-additive on some $\sigma $-algebra $\sf A$
of subsets of $X$. Moreover, a correlation operator $B$ is of
class $L_1$, that is, $Tr (B)<\infty $, if and only if
each finite dimensional over $\bf K$
projection of $\mu $ is a $\sigma $-additive $q$-Gaussian Borel measure.}
\par {\bf Proof.} From Definition 2.13 it follows, that
each one dimensional over $\bf K$ projection $\mu _{x\bf K}$
of a measure $\mu $ is $\sigma $-additive on the Borel
$\sigma $-algebra $Bf ({\bf K})$, where $0\ne x=e_{k,l}\in X_l$.
Therefore, $\mu $ is defined and finite additive on a cylindrical
algebra ${\sf U}:=\bigcup_{k_1,...,k_n;l}\phi _l^{-1}
[(\phi ^l_{k_1,...,k_n})^{-1}
(Bf (sp_{\bf K}\{ e_{k_1,l},...,e_{k_n,l} \} ))]$, where
$\phi ^l_{k_1,...,k_n}: X_l\to sp_{\bf K}(e_{k_1,l},...,e_{k_n,l})$
is a projection on a $\bf K$-linear span of vectors $e_{k_1,l},...,
e_{k_n,l}$. This means that $\mu $ is a bounded quasimeasure
on $\sf U$. Since $\hat \mu $ is the positive definite
function, then $\mu $ is realvalued. In view of the non-Archimedean
analog of the Bochner-Kolmogorov theorem (see I.2.27 \cite{lulapm})
$\mu $ has an extension to a $\sigma $-additive probability measure on
a $\sigma $-algbera $\sigma \sf U$, that is, a minimal $\sigma $-algebra
of subsets of $X$ containing $\sf U$.
If $J: X_j\to X_j$ is a $\bf K$-linear operator diagonal
in the basis $\{ e_{k,j}: k \} $, then
for $z$ such that $z(x)=z_j(\phi _j(x))$
for each $x\in X$ and a symmetric nonegative definite operator $F$
as in \S 2.13
\par $(i)$ $F(v_q(z\circ J),v_q(z\circ J))=E(v_q(z),v_q(z))$, where
\par $(ii)$ $E_{k,l}=F_{k,l}|J_{k,k}|^{q/2}|J_{l,l}|^{q/2}$ for each
$k, l\in \alpha _j$.
If $F\in L_a$ (that is, $F^a\in L_1$) and $J\in L_q$
(that is, $diag (v_1(J_{l,l}): l)\in L_q$), then
\par $(iii)$ $E\in L_{aq/(a+q)}$ for each $a>0$
(see Theorem 8.2.7 \cite{pietsch}).
In particular, taking $a$ tending to $\infty $ and
$F=I$ we get $E\in L_q$, since $L_{\infty }$ is the space
of bounded linear operators.
Using the orthonormal bases in $X_j$ for each $j$
we get the embedding of $X_j$ into its topologically adjoint space
$X_j^*$ of all continuous $\bf K$-linear functionals on $X_j$.
For each $z\in X^*$ there exists a non-Archimedean direct
sum decomposition $X=X_z\oplus ker (z)$, where
$X_z$ is a one dimensional over $\bf K$ subspace in $X$.
Therefore, the set ${\sf D}_{q,B,X}$ of functionals $z$ on $X$ from
\S 2.13 separates points of $X$.
If for a given one dimensional over $\bf K$ subspace $W$
in $X$ it is the equality $B(v_q(z),v_q(z))=0$ for each
$z\in W$, then the projection $\mu _W$ of $\mu $
is the atomic measure with one atom being a singleton.
If $B\in L_1$, then $B(v_q(z),v_q(z))$ and hence $\hat \mu (z)$ is correctly
defined for each $z\in {\sf D}_{q,B,X}$.
From Equalities $(i,ii)$, Inclusion $(iii)$ above and
analogously to Theorem II.2.2 \cite{lustpr} we get
the second statement of this theorem.
\par {\bf 2.15. Corollary.} {\it Let $X$ be a complete locally
$\bf K$-convex space of separable type over a local field $\bf K$,
then for each constant $q>0$ there exists a nondegenerate
symmetric positive definite operator $B\in L_1$ such that
a $q$-Gaussian measure is $\sigma $-additive on $Bf (X)$ and
each its one dimensional over $\bf K$ projection is absolutely
continuous relative to the nonnegative Haar measure on $\bf K$.}
\par {\bf Proof.} A space $Y$ from \S 2.12 corresponding to
$X$ is a separable locally $\bf R$-convex space.
Therefore, $Y$ in a weak topology is isomorphic with
${\bf R}^{\aleph _0}$ from which the existence of $B$ follows.
For each $\bf K$-linear finite dimensional over $\bf K$ subspace
$S$ a projection $\mu ^S$ of $\mu $ on $S\subset X$
exists and its density $\mu ^S(dx)/w(dx)$ relative to the
nonegative nondegenerate Haar measure $w$ on $S$ is the inverse
Fourier transform $F^{-1}({\hat \mu }|_{S^*})$
of the restriction of ${\hat \mu }$ on $S^*$ (see about
the Fourier transform on non-Archimedean spaces \S VII \cite{vla3}).
For such $B$ each one dimensional projection of $\mu $
corresponding to $\hat \mu $ has a density that is a
continuous function belonging to $L^1({\bf K},w,Bf ({\bf K}),{\bf R})$,
where $w$ denotes the nonnegative Haar measure on $\bf K$.
\par {\bf 2.16. Proposition.} {\it Let $\mu _{q,B,\gamma }$ and
$\mu _{q,E,\delta }$ be two $q$-Gaussian measures with
correlation operators $B$ and
$E$ of class $L_1$, then there exists a convolution
of these measures $\mu _{q,B,\gamma }*\mu _{q,E,\delta }$, which is
a $q$-Gaussian measure $\mu _{q,B+E,\gamma +\delta }$.}
\par {\bf Proof.} Since $B$ and $E$ are nonnegative, then
$(B+E)(y,y)=B(y,y)+E(y,y)\ge 0$ for each $y\in Y$, that is,
$B+E$ is nonegative. Evidently, $B+E$ is symmetric.
In view of \cite{pietsch} $B+E$ is of class $L_1$.
Therefore, $\mu _{q,B+E,\gamma +\delta }$ is the $\sigma $-additive
$q$-Gaussian measure together with $\mu _{q,B,\gamma }$ and
$\mu _{q,E,\delta }$ in accordance with Theorem 2.14.
Moreover, $\mu _{q,B+E,\gamma +\delta }$ is defined on the $\sigma $-algebra
$\sigma {\sf U}_{B+E}$ containing the union of $\sigma $-algebras
$\sigma {\sf U}_B$ and $\sigma {\sf U}_E$
on which $\mu _{q,B,\gamma }$ and $\mu _{q,E,\delta }$
are defined correspondingly,
since $ker (B+E)\subset ker (B)\cap ker (E)$.
Since ${\hat \mu }_{q,B+E,\gamma +\delta }={\hat \mu }_{q,B,\gamma }
{\hat \mu }_{q,E,\delta }$,
then $\mu _{q,B+E,\gamma +\delta }=\mu _{q,B,\gamma }*\mu _{q,E,\delta }$
(see Proposition I.2.11 \cite{lulapm} and use projective limits).
\par {\bf 2.17. Definitions.} Let $B$ and $q$ be as in
\S 2.14 and denote by $\mu _{q,B,\gamma }$ the corresponding $q$-Gaussian
measure on $H$. Let $\xi $ be a stochastic process with a real time
$t\in T\subset \bf R$ (see \S 2.2),
then it is called a non-Archimedean $q$-Wiener process
with real time, if
\par $(ii)'$ the random variable $\xi (t,\omega )-\xi (u,\omega )$ has
a distribution $\mu _{q,(t-u)B,\gamma }$ for each $t\ne u\in T$. 
\par Let $\xi $ be a stochastic process with a non-Archimedean time
$t\in T\subset \bf F$, where $\bf F$ is a local field,
then $\xi $ is called a non-Archimedean $q$-Wiener process
with $\bf F$-time, if
\par $(ii)"$ the random variable $\xi (t,\omega )-\xi (u,\omega )$ has
a distribution $\mu _{q,ln [\chi _{\bf F}(t-u)]B,\gamma }$
for each $t\ne u\in T$,
where $\chi _{\bf F}: {\bf F}\to S^1$ is a continuous character
of $\bf F$ as the additive group.
\par {\bf 2.18. Proposition.} {\it For each given $q$-Gaussian measure
a non-Archimedean $q$-Wiener process with real ($\bf F$ respectively)
time exists.}
\par {\bf Proof.} In view of Proposition 2.16
for each $t>u>b$ a random variable $\xi (t,\omega )-\xi (b,\omega )$
has a distribution $\mu _{q,(t-b)B,\gamma }$ for real time parameter.
If $t$, $u$, $b$ are pairwise different points in $\bf F$,
then $\xi (t,\omega )-\xi (b,\omega )$ has a distribution
$\mu _{q,ln [\chi _{\bf F}(t-b)]B,\gamma }$, since $ln [\chi _{\bf F}(t-u)]+
ln [\chi _{\bf F}(u-b)]=ln [\chi _{\bf F}(t-b)]$.
This induces the Markov quasimeasure $\mu ^{(q)}_{x_0,\tau }$
on $(\prod_{t\in T}(H_t,{\sf U}_t)),$ where
$H_t=H$ and ${\sf U}_t=Bf (H)$ for each $t\in T$
(see \S VI.1.1 \cite{dalfo} and \S 3 in I \cite{lustpr}).
Therefore, the Chapman-Kolmogorov equation is accomplished:
$$P(b,x,t,A)=\int_H P(b,x,u,dy) P(u,y,t,A)$$
for each $A\in Bf (H)$. An abstract probability space
$(\Omega ,{\sf F},\lambda )$ exists due to the Kolmogorov theorem,
hence the corresponding space $L^r$ exists.
Therefore, conditions of defintions 2.2 and 2.17 are satisfied.
\par {\bf 2.19. Proposition.} {\it Let $\xi $ be a $q$-Gaussian process
with values in a Banach space $H=c_0(\alpha ,{\bf K})$
a time parameter $t\in T$ and a positive definite
correlation operator $B$ of trace class
and $\gamma =0$, where $card (\alpha )\le \aleph _0$,
either $T\subset \bf R$ or $T\subset \bf F$. Then either
$$(i)\quad \lim_{N\in \alpha }M_t \| (v_q(e^1(\xi (t,\omega )),...,
v_q(e^N(\xi (t,\omega ))) \|_{l_2}^2=t Tr(B) \mbox{ or }$$
$$(ii)\quad \lim_{N\in \alpha }M_t \| (v_q(e^1(\xi (t,\omega )),...,
v_q(e^N(\xi (t,\omega )) \|_{l_2}^2=[ln (\chi _{\bf F}(t))] Tr(B)
\mbox{ respectively}.$$ }
\par {\bf Proof.} At first we consider moments of a $q$-Gaussian measure
$\mu _{q,B,\gamma }$. We define moments
$m^q_k(e^{j_1},...,e^{j_k}):=
\int_Hv_{2q}(e^{j_1}(x))...v_{2q}(e^{j_k}(x))\mu _{q,B,\gamma }(dx)$
for linear continuous functionals $e^{j_1},...,e^{j_k}$
on $H$ such that $e^l(e_j)=\delta ^l_j$, where in our previous
notation $\{ e_j: j\in \alpha \} $ is the standard orthonormal base in $H$.
\par Consider partial pseudodifferential operators
$\mbox{ }_P\partial ^u_j$ given by the equation
$$(iii)\quad \mbox{ }_P\partial ^u_j\psi (x):={\sf F}_j^{-1}(
|\tilde x_j|_{\bf K}^u {\hat \psi }(\tilde x))(x),$$ where
the norm $|b|_{\bf K}=mod_{\bf K}(b)$ on $\bf K$ is chosen coinciding
with the modular function associated with the nonnegative
nondegenerate Haar measure $w$ on $\bf K$ (about the modular function
see \cite{wei}), $u \in {\bf C}\setminus \{ -1 \} $,
${\hat \psi }:={\sf F}_j(\psi )$ is the Fourier transform of $\psi $
by a variable $x_j\in \bf K$ such that ${\sf F}_j$ is defined
relative to the Haar measure $w$ on $\bf K$ \cite{vla3}.
From the change of variables formula
$\int_{\bf K}f(ax+b)g(x)w(dx)=\int_{\bf K}f(y)g((y-b)/a)|a|_{\bf K}^{-1}
w(dy)$ for each $f$ and $g\in L^2({\bf K},Bf ({\bf K}),w,{\bf C})$,
$a\ne 0$ and $b\in \bf K$, also the Fubini theorem and the Fourier
transform on $\bf K$ it follows that $f_{-\alpha }*f_{u+1}=f_{u+1-\alpha }$
for $u\ne \alpha $ and $\Gamma _{\bf K}(u+1)
|\xi _j|_{\bf K}^{-u-1}={\sf F}(|x_j|^u)$, where
$f_u(x_j):=|x_j|^{u-1}/\Gamma _{\bf K}(u)$,
$\Gamma _{\bf K}$ is the non-Archimedean gamma function,
$\Gamma _{\bf K}(u):=\int_{\bf K}|z|_{\bf K}^{u-1}\chi (z)w(dz)$,
$\chi : {\bf K}\to S^1$ is the character of $\bf K$ as the additive group
such that $\chi (z):=\prod_{j=1}^m\chi _p({z'}_j)$,
${z'}_j\in \bf Q_p$, $z=({z'}_1,...,{z'}_m)\in \bf K$ for $\bf K$
considered as the $\bf Q_p$-linear space, $m\in \bf N$,
$dim_{\bf Q_p}{\bf K}=m$, $\chi _p: {\bf Q_p}\to S^1$
is the standard character such that $\chi _p(y):=\exp (2\pi
i \{ y \} _p)$, $\{ y \} _p:=\sum_{l<0}a_lp^l$ for $|y|_{\bf Q_p}>1$
and $\{ y \} _p=0$ for $|y|_{\bf Q_p}\le 1$,
$y=\sum_la_lp^l$, $a_l\in \{ 0,1,...,p-1 \} $, $l\in \bf Z$,
$\min (l: a_l\ne 0 )=:ord_p(y)>-\infty $.
Therefore, $\mbox{ }_P\partial ^u_j|x_j|^n=|x_j|^{n-u}
\Gamma _{\bf K}(n)/\Gamma _{\bf K}(n-u)$, where $n\in {\bf C}
\setminus \{ -1 \} $.
A function $\psi $ for which $\mbox{ }_P\partial ^u_j\psi $ exists
is called pseudodifferentiable of order $u$ by variable $x_j$.
\par From $m^{q/2}_k(e^{j_1},...,e^{j_k})={\sf F}^{-1}
(|x_{j_1}|^{q/2}...|x_{j_k}|^{q/2}{\sf F}(\mu ))(0)$,
since ${\sf F}(hg)={\sf F}(h)*{\sf F}(g)$ for functions $h$ and $g$
in the Hilbert space $L^2({\bf K},Bf ({\bf K}),w,{\bf C})$
it follows that $m^{q/2}_{2k}(e^{j_1},...,e^{j_{2k}})=
\mbox{ }_P\partial ^{q/2}_{j_1}...\mbox{ }_P\partial ^{q/2}_{j_{2k}}
{\hat \mu }(0)=([\mbox{ }_PD^{q/2}]^{2k}{\hat \mu }(0)).(e_{j_1},...,
e_{j_{2k}})$, where $\mbox{ }_PD^{q/2}$
is a $\bf K$-linear pseudodifferential
operator by $x\in H$ such that $(\mbox{ }_PD^{q/2}\psi (x)).e_j:=
\mbox{ }_P\partial ^{q/2}_j\psi (x)$.
Then
\par $(iv)$ $m^{q/2}_{2n}(e^{j_1},...,e^{j_{2n}})=$
$(-1)^n(n!)^{-1}2^{-n} [\mbox{ }_PD^{q/2}]^{2n}[B(v_q(z),v_q(z)]^n.
(e_{j_1},...,e_{j_{2n}})$
\par $=(n!)^{-1}2^{-n}\sum_{\sigma \in \Sigma _{2n}}
B_{\sigma (j_1),\sigma (j_2)}...B_{\sigma (j_{2n-1}),\sigma (j_{2n})}$, \\
since $\gamma =0$ and $\chi _{\gamma }(z)=1$,
where $\Sigma _k$ is the symmetric group of
all bijective mappings $\sigma $ of the set $\{ 1,...,k \} $
onto itself, $B_{l,j}:=B(e_j,e_l)$, since $Y^*=Y$ for
$Y=l_2(\alpha ,{\bf R})$.
Therefore, for each $B\in L_1$ and $A\in L_{\infty }$
we have $\int_HA(v_q(x),v_q(x))\mu _{q,B,0}(dx)=$
$\lim_{N\in \alpha }\sum_{j=1}^N\sum_{k=1}^N
A_{j,k}m^{q/2}_2(e_j,e_k)$ $=Tr(AB)$.
\par In particular for $A=I$ and $\mu _{q,tB,0}$ corresponding
to the transition measure of $\xi (t,\omega )$
we get Formula $(i)$ for a real time parameter,
using $\mu _{q,ln [\chi _{\bf F}(t)]B,0}$
we get Formula $(ii)$ for a time parameter belonging to $\bf F$,
since $\xi (t_0,\omega )=0$ for each $\omega $.
\par {\bf 2.20. Corollary.} {\it Let $H=\bf K$
and $\xi $, $B=1$, $\gamma $ be as in Proposition 2.19, then
$$(i)\quad M(\int_{t\in [a,b]}
\phi (t,\omega )|d\xi (t,\omega )|_{\bf K}^q)=
M[\int_a^b\phi (t,\omega )dt]$$
for each $a<b\in T$ with real time, where
$\phi (t,\omega )\in L^2(\Omega ,{\sf U},\lambda ,C^0_0(T,{\bf R}))$
$\xi \in L^q(\Omega ,{\sf U},\lambda ,C^0_0(T,{\bf K}))$,
$(\Omega ,{\sf U},\lambda )$ is a probability measure space.}
\par {\bf Proof.} Since $\int_{t\in [a,b]}\phi (t,\omega )
|d\xi (t,\omega )|_{\bf K}^q$ \\
$=\lim_{\max_j (t_{j+1}-t_j)\to 0}
\sum_{j=1}^N\phi (t_j,\omega ) |\xi (t_{j+1},\omega )-
\xi (t_j,\omega )|_{\bf K}^q$ for $\lambda $-almost all $\omega
\in \Omega $, then applying Fromula $2.19.(i)$ to each
$|\xi (t_{j+1},\omega )-\xi (t_j,\omega )|_{\bf K}^q$
and taking the limit by finite partitions $a=t_1<t_2<...<t_{N+1}=b$
of the segment $[a,b]$ we get Formula $2.20.(i)$.
\par {\bf 2.21. Remarks.} In the classical case with $q=2$
and $\bf R$ instead of $\bf K$ there is analogous formula
$M([\int_{t\in [a,b]}\phi (t,\omega )dB_t(\omega )]^2)=
M[\int_a^b\phi (t,\omega )^2dt$ known as the It$\hat o$ formula
(see the classical case in \cite{boui,coxmil,dalfo,gihsko,
gisk3,ikwat,itkean,malb,oksb}).
Another analogs of the It$\hat o$ formula were given in
II \cite{lustpr}. Certainly it is impossible to get
in the non-Archimedean case all the same properties of Gaussian measures
and Wiener process (Brownian motion) as in the classical case.
Therefore, there are different possibilities for seeking
non-Archimedean analogs of Gaussian measures and Wiener processes
depending on a set of properties supplied with these objects.
Giving our definitions we had the intention to take into
account the most important properties.
\par Since ${\sf F}(\chi _{\gamma })(y)=\delta (y-\gamma )$
and $[\delta (y-\gamma )*h(y)](x)=h(x-\gamma )$
for any continuous function $h$,
then $\int_H|x_{j_1}-\gamma _{j_1}|^{q/2}...
|x_{j_k}-\gamma _{j_k}|^{q/2}d\mu _{q,B,\gamma }=$
$\int_H|x_{j_1}|^{q/2}...|x_{j_k}|^{q/2}d\mu _{q,B,0}$,
consequently, $\gamma $ plays in some sence the mean value role.
\par If $A>0$ on $Y=l_2(\alpha ,{\bf K})$, then
\par $\mu _{q,B,0} \{ x: A(v_q(x),v_q(x))\ge 1 \} \le Tr (AB)$ and
\par $\mu _{q,B,0} \{ x: |A(v_q(x),v_q(x))-Tr (AB)|\le c (Tr (AB))^{1/2}
\} \ge 1-2 \| AB\| /c^2$ for each $c>0$ due to the Chebyshev
inequality and Formula $2.19.(iv)$.
\par {\bf 2.22. Definitions and Notes.} Consider a pseudodifferential operator
on $H=c_0(\alpha ,{\bf K})$ such that
$$(i)\quad {\sf A}=\sum_{0\le k\in {\bf Z};
j_1,...,j_k\in \alpha }(-i)^kb^k_{j_1,...,j_k}
\mbox{ }_P\partial _{j_1}...\mbox{ }_P\partial _{j_k},$$
where $b^k_{j_1,...,j_k}\in \bf R$,
$\mbox{ }_P\partial _{j_k}:=\mbox{ }_P\partial ^1_{j_k}$.
If there exists $n:=\max \{ k: b^k_{j_1,...,j_k}
\ne 0, j_1,...,j_k\in \alpha  \} $,
then $n$ is called an order of ${\sf A}$, $Ord ({\sf A})$.
If ${\sf A}=0$, then by definition $Ord ({\sf A})=0$.
If there is not any such finite $n$, then $Ord ({\sf A})=\infty $.
We suppose that the corresponding form $\tilde A$ on $\bigoplus_kY^k$
is continuous into $\bf C$, where
$$(ii)\quad {\tilde A}(y)=-\sum_{0\le k\in {\bf Z};
j_1,...,j_k\in \alpha }(-i)^kb^k_{j_1,...,j_k}y_{j_1}...y_{j_k},$$
$y\in l_2(\alpha ,{\bf R})=:Y$.
If ${\tilde A}(y)>0$ for each $y\ne 0$ in $Y$, then
$\sf A$ is called strictly elliptic pseudodifferential operator.
The phase multiplier $(-i)^k$ is inserted into
the definition of $\sf A$ for in the definition of $\mbox{ }_P\partial _j$
it was omitted in comparison with the classical case.
\par Let $X$ be a complete locally $\bf K$-convex space,
let $Y$ be a corresponding complete locally $\bf R$-convex space
(see \S 2.12), let $Z$ be a complete locally $\bf C$-convex space.
For $0\le n\in \bf R$ a space of all functions $f: X\to Z$
such that $f(x)$ and $(\mbox{ }_PD^kf(x)).(y^1,...,y^{l(k)})$ are continuous
functions on $X\times Y^{l(k)}$, $l(k):=[k]+sign \{ k \} $ for each $k\in
\bf N$ such that $k\le [n]$ and also for $k=n$
is denoted by $\mbox{ }_P{\cal C}^n(X,Y,Z)$ and $f\in \mbox{ }_P
{\cal C}^n(X,Y,Z)$
is called $n$ times continuously pseudodifferentiable, where
$[n]\le n$ is an integer part of $n$, $1> \{ n \} :=n-[n]\ge 0$
is a fractional part of $n$, $sign (b)=1$ for each $b>0$, $sign (0)=0$,
$sign (b)=-1$ for $b<0$, $y^1,...,y^{l(k)}\in Y$.
Then $\mbox{ }_P{\cal C}^{\infty }(X,Y,Z):=\bigcap_{n=1}^{\infty }
\mbox{ }_P{\cal C}^n(X,Y,Z)$
denotes a space of all infinitely pseudodifferentiable functions.
\par {\bf 2.23. Theorem.} {\it Let $\sf A$ be a
strictly elliptic pseudodifferential
operator on $H=c_0(\alpha ,{\bf K})$, $card (\alpha )\le \aleph _0$,
and let $t\in T=[0,b]\subset \bf R$. Suppose also that
$u_0(x-y)\in L^2(H,Bf(H),\mu _{t{\tilde A}},{\bf C})$
for each marked $y\in H$ as a function by $x\in H$,
$u_0(x)\in \mbox{ }_P{\cal C}^{Ord ({\sf A})}(H,Y,{\bf C})$.
Then the non-Archimedean analog of the Cauchy problem
$$(i)\quad \partial u(t,x)/\partial t={\sf A}u,\quad u(0,x)=u_0(x)$$
has a solution given by
$$(ii)\quad u(t,x)=\int_Hu_0(x-y)\mu _{t{\tilde A}}(dy),$$
where $\mu _{t{\tilde A}}$ is a $\sigma $-additive Borel measure on $H$
with a characteristic functional ${\hat \mu }_{t{\tilde A}}(z)
:=\exp [-t{\tilde A}(v_1(z))]$.}
\par {\bf Proof.} In accordance with \S \S 2.12 and 2.22
we have $Y=l_2(\alpha ,{\bf R})$.
In view of the conditions of this theorem
the function $\exp [-t{\tilde A}(v_1(z))]$
is continuous on $H\hookrightarrow H^*$ for each $t\in \bf R$
such that the family $H$ of continuous $\bf K$-linear functionals
on $H$ separates points in $H$. In view of the Minlos-Sazonov theorem
I.2.35 \cite{lulapm} it defines a $\sigma $-additive Borel measure
on $H$ for each $t>0$ and hence for each $t\in (0,b]$.
The functional $\tilde A$ on each ball of radius $0<R<\infty $
in $Y$ is a uniform limit of its restrictions
${\tilde A}|_{\bigoplus_k[sp_{\bf K}(e_1,...,e_n)]^k},$
when $n$ tends to the infinity, since $\tilde A$ is
continuous on $\bigoplus_kY^k$.
Since $u_0(x-y)\in L^2(H,Bf(H),\mu _{t{\tilde A}},{\bf C})$
and a space of cylindrical functions is dense in the latter Hilbert space,
then due to the Parceval-Steclov equality and the Fubini theorem
it follows that $\lim_{P\to I}{\sf F}_{Px}u_0(Px)){\hat \mu }_{t\tilde A}
(y+Px)$ converges in $L^2(H,Bf(H),\mu _{t{\tilde A}},{\bf C})$
for each $t$, since $\mu _{t_1\tilde A}*\mu _{t_2\tilde A}=\mu_{(t_1+t_2)
\tilde A}$ for each $t_1$, $t_2$ and $t_1+t_2\in T$, where
$P$ is a projection on a finite dimensional over $\bf K$ subspace
$H_P:=P(H)$ in $H$, $H_P\hookrightarrow H$,
$P$ tends to the unit operator $I$ in the strong operator
topology, ${\sf F}_{Px}u_0(Px)$ denotes a Fourier transform by
the variable $Px\in H_P$.
Consider a function $v:={\sf F}_x(u)$, then $\partial v(t,x)/\partial t
=-{\tilde A}(v_1(x))v(t,x)$, consequently,
$v(t,x)=v_0(x)\exp [-t{\tilde A}(v_1(x))]$. From $u(t,x)={\sf F}_x^{-1}
(v(t,x))$, where as above ${\sf F}_x(u)$ denotes the Fourier
transform by the variable $x\in H$ such that
${\sf F}_x(u(t,x))=\lim_{n\to \infty }{\sf F}_{x_1,...,x_n}u(t,x)$.
Therefore, $u(t,x)=u_0(x)*{\sf F}_x^{-1}({\hat \mu }_{t{\tilde A}})
=$ $\int_Hu_0(x-y)\mu _{t{\tilde A}}(dy)$,
since $u_0(x-y)\in L^2(H,Bf(H),\mu _{t{\tilde A}},{\bf C})$
and $\mu _{t{\tilde A}}$ is the bounded measure
on $Bf(H)$ and $|\int_Hu_0(x-y)\mu _{t{\tilde A}}(dy)|\le $
$(\int_H |u_0(x-y)|^2\mu _{t{\tilde A}}(dy))\mu _{t{\tilde A}}(H)
<\infty $.
\par {\bf 2.24. Note.} In the particular case of $Ord ({\sf A})=2$
and ${\tilde A}$ corresponding to the Laplace operator, that is,
${\tilde A}(y)=\sum_{l,j} g_{l,j}y_ly_j$, Equation
$2.23.(i)$ is (the non-Archimedean analog of) the heat equation
on $H$. This provides the interpretation of the $2$-Gaussian
measure $\mu _{t\tilde A}=\mu _{2,t{\tilde A},0}$.
For $dim_{\bf K}H<\infty $ the density $\mu _{t\tilde A}(dx)/w(dx)$
is called the heat kernel, where $w$ is the nonnegative nondegenerate
Haar measure on $H$.
\par For $Ord ({\sf A})<\infty $ the form ${\tilde A}_0(y)$ corresponding
to sum of terms with $k=Ord ({\sf A})$ in Formula $2.22.(ii)$ is called
the principal symbol of operator $\sf A$. If ${\tilde A}_0(y)>0$
for each $y\ne 0$, then $\sf A$ is called an elliptic pseudodifferential
operator. Evidently, Theorem 2.23 is true for elliptic $\sf A$
of $Ord ({\sf A})<\infty $, since $\exp [-t {\tilde A}(v_1(z))]$
is the bounded continuous realvalued positive definite function.
\par {\bf 2.25. Remark and Definitions.} Let linear spaces $X$ over $\bf K$
and $Y$ over $\bf R$ be as in \S 2.12 and $B$ be a symmetric
nonegative definite (bilinear) operator on a dense $\bf R$-linear subspace
${\sf D}_{B,Y}$ in $Y^*$. A quasimeasure $\mu $ with a characteristic
functional
$${\hat \mu }(\zeta ,x):=\exp [-\zeta B(v_q(z),v_q(z))/2]\chi _{\gamma }(z)$$
for a parameter $\zeta \in \bf C$ with $Re (\zeta )\ge 0$
defined on ${\sf D}_{q,B,X}$
is called a complexvalued Gaussian measure and is denoted
by $\mu _{q,\zeta B,\gamma }$ also, where ${\sf D}_{q,B,X}:=
\{ z\in X^*:$ $\mbox{there exists}$ $j\in \Upsilon $
$\mbox{such that}$ $z(x)=z_j(\phi _j(x))$ $\forall x\in X,$
$ v_q(z)\in {\sf D}_{B,Y} \} $.
\par {\bf 2.26. Proposition.} {\it Let $X={\sf D}_{q,B,X}$
and $B$ be positive definite, then
for each function $f(z):=\int_X\chi _z(x)\nu (dx)$ with
a complexvalued measure $\nu $ of finite variation
and each $Re (\zeta )>0$ there exists
$$(i)\quad \int_Xf(z)\mu _{\zeta B}(dz)=
\lim_{P\to I}\int_Xf(Pz)\mu ^{(P)}_{\zeta B}(dz)$$
$$=\int_X\exp (-\zeta B(v_q(z),v_q(z))/2)\chi _{\gamma }(z)\nu (dz),$$
where $\mu ^{(P)}(P^{-1}(A)):=\mu (P^{-1}(A))$ for each
$A\in Bf (X_P)$, $P: X\to X_P$ is a projection on a
$\bf K$-linear subspace $X_P$, a convergence $P\to I$ is
considered relative to a strong operator topology.}
\par {\bf Proof.} A complexvalued measure $\nu $ can be presented
as $\nu =\nu _1-\nu _2+i\nu _3 -i\nu _4$, where
$\nu _j$ are nonnegative measures, $j=1,...,4$, $i=(-1)^{1/2}$.
Using the projective limit decomposition of $X$ and \S 2.27
in I \cite{lulapm} we get that \\
$(ii)\quad \int_Xf(z)\mu _{\zeta B}(dz)=
\lim_{P\to I}\int_Xf(Pz)\mu ^{(P)}_{\zeta B}(dz)$.
On the other hand, for each
finite dimensional over $\bf K$ subspace $X_P$ \\
$(iii)\quad \int_Xf(Pz)\mu ^{(P)}_{\zeta B}(dz)=\int_{X_P}
\{ \exp (-\zeta B(v_q(z),v_q(z))/2) \chi _{\gamma }(z) \}
|_{X_P}\nu ^{X_P}(dz).$
Since each measure $\nu _j$ is nonegative and finite,
then due to Lemma 2.3 and \S 2.6 in I \cite{lulapm} there exists the limit \\
$\lim_{P\to I} \int_{X_P}\{ \exp [-\zeta B(v_q(z),v_q(z))/2]
\chi _{\gamma }(z) \} |_{X_P}
\nu ^{X_P}(dz)$\\
$= \int_X\exp (-\zeta B(v_q(z),v_q(z))/2)\chi _{\gamma }(z)\nu (dz).$
\par {\bf 2.27. Proposition.} {\it If conditions of Proposition
2.26 are satisfied and
$$(i)\quad \int_{X_P} |f(Px)|w^{X_P}(dx)<\infty $$
for each finite dimensional over $\bf K$ subspace $X_P$ in $X$,
then Formula $2.26.(i)$ is accomplished for $\zeta $ with
$Re (\zeta )=0$, where $w^{X_P}$ is a nonegative nondegenerate
Haar measure on $X_P$.}
\par {\bf Proof.} The finite dimensional over $\bf K$ distribution
$\mu ^{X_P}_{q,iB,\gamma }/w^{X_P}(dx)={\sf F}^{-1}(
{\hat \mu }_{q,iB,\gamma }|_{X_P})$ is locally $w^{X_P}$-integrable, but
does not belong to the space $L^1(X_P,Bf (X_P),w^{X_P},{\bf C})$.
In view of Condition $2.27.(i)$ above
and the Fubini theorem and using the Fourier transform of generalized
functions (see \S VII.3 \cite{vla3}) we get
Formulas $2.26.(ii,iii)$. Taking the limit by $P\to I$
we get Formula $2.26.(i)$ in the sence of distributions.
\par {\bf 2.28. Remark.} A measure $\mu _{q,iB,\gamma }$ is
the non-Archimedean analog of the Feynman quasimeasure. Put
$$(i)\quad \mbox{ }_F\int_X
f(x)\mu _{q,iB,\gamma }(dx):=\lim_{\zeta \to i}
\int_Xf(x)\mu _{q,\zeta B,\gamma }(dx)$$
if such limit exists. If conditions of Proposition 2.26 are
satisfied, then $\psi (\zeta ):=\int_Xf(x)\mu _{q,\zeta B,\gamma }(dx)$
is the holomorphic function on $\{ \zeta \in {\bf C}: Re (\zeta )>0 \} $
and it is continuous on $\{ \zeta \in {\bf C}: Re (\zeta )\ge 0 \} $,
consequently,
$$(ii)\quad \mbox{ }_F\int_X f(x)\mu_{q,iB,\gamma }(dx)=
\int_X\exp \{ -i B(v_q(x),v_q(x))/2 \} \chi _{\gamma }(x)\nu (dx).$$
\section{Stochastic processes on non-Archimedean manifolds.}
To avoid misunderstanding at first definitions and notations are given.
\par {\bf 3.1. Definitions and Notes.} Let $M$ be a $C^n$-manifold
on a Banach space $X$ over a non-Archimedean field $\bf K$
complete relative to its norm
with an atlas $At (M):= \{ (U_j,\phi _j): j\in \Lambda _M \} $
such that $U_j$ is an open covering of $M$ and
$\phi _j: U_j \to \phi _j(U_j)$ is a homeomorphism,
$\phi _j(U_j)$ is open in $X$, $\phi _l\circ \phi _j^{-1}:
\phi _j(U_j\cap U_l)\to \phi _l(U_j\cap U_l)$ is a diffeomorphism
of class $C^n$ for each $U_j\cap U_l\ne \emptyset $,
the space $C^n(U,Y)$ of functions from an open subset $U\subset X$
into a Banach space $Y$ over $\bf K$ is defined in terms of
difference quotients (see \cite{sch1,lubp2}).
\par Since the derivative $(\phi _l\circ
\phi _j^{-1})'(x)$ is a linear continuous operator
on $\phi _j(U_j\cap U_l)\times X$ of class $C^{n-1}$
for each $n\ge 1$ and there exists a derivative of an inverse
operator $(\phi _j\circ \phi _l^{-1})'(y)$ on
$\phi _l(U_j\cap U_l)\times X$, then $(\phi _l\circ \phi _j^{-1})'(x)
\in GL(X)$ for each $x\in \phi _j(U_j\cap U_l)$,
where $GL(X)$ is the group of invertable $\bf K$-linear bounded operators
of $X$ onto $X$.
Therefore, for each $n\ge 1$ there exists a functor $T$
such that $T(\phi _l\circ \phi _j^{-1})(x):=
(\phi _{l,j}(x),\phi _{l,j}'(x))$ for each $x\in
(\phi _l\circ \phi _j^{-1})(U_j\cap U_l)$,
$T(\phi _j(U_j)):=\phi _j(U_j)\times X$,
where $\phi _{l,j}:=(\phi _l\circ \phi _j^{-1})$.
For $n\ge 1$ put $TM=\bigcup_{j\in \Lambda _M}TU_j$ with the atlas
$At (TM):= \{ (U_j\times X,T\phi _j): j\in \Lambda _M \} $
such that $T\phi _j: TU_j \to \phi _j(U_j)\times X$
is a homeomorphism, $T\phi _j|_{\{ x \} \times X}=:T_x\phi _j$
is a bounded continuous operator on $X$ by the second argument for each
$x\in U_j$. Thus $TM=\bigcup_{x\in M}T_xM$, where $T_x\phi _j:
T_xU_j\to T_{\phi _j(x)}\phi _j(U_j)$ is a $\bf K$-linear
isomorphism for each $j\in \Upsilon _M$, where
$T_{\phi _j(x)}\phi _j(U_j)=\{ \phi _j(x) \} \times X$.
$TM$ is called the total tangent space of $M$,
$T_xM$ is called the tangent space of $M$ at $x$.
The projection $\tau :=\tau _M: TM\to M$ is given by
$\tau _M(s)=x$ for each vector $s\in T_xM$, $\tau _M$ is called
the tangent bundle.
\par {\bf 3.1.1} If $M$ and $N$
are two $C^l$-manifolds on Banach spaces $X$ and $S$
over $\bf K$ with $l\ge n$,
where $At (N):= \{ (V_j,\psi _j): j\in \Lambda _N \} $
and $f: M\to N$ is a continuous mapping, then by the definition
$f\in C^n(M,N)$, if $\psi _l\circ f\circ \phi _j^{-1}\in C^n(W_{l,j},Y)$
for each $W_{l,j}:=\phi _j(f^{-1}(V_l)\cap U_j)\ne \emptyset $.
A norm in $C^n(X,S)$ induces a complete uniformity in $C^n(M,N)$.
If $n\ge 1$ and $f\in C^n(M,N)$, then there exists
$Tf: TM\to TN$ and $Tf\in C^{n-1}(TM,TN)$.
\par {\bf 3.1.2.} Let $H$ and $X$ be two Banach spaces over
a non-Archimedean field $\bf K$. Let $M$ be a $C^l$-manifold
on $X$ and let $P$ be a manifold with a mapping $\pi : P\to M$
such that $\pi $ is surjective and $\pi ^{-1}(x)=:P_x=:H_x$ is
a Banach space over $\bf K$ isomorphic to $H$ for each $x\in M$,
$\pi $ is called a projection, $\pi ^{-1}(x)$ is called a fibre
of $\pi $ over $x$. Suppose that $P$ is supplied with an atlas
$At(P)= \{ (U_j,\phi _j,P\phi _j): j\in \Lambda _M \} $
consistent with $At (M)$ such that $pr_1\circ P\phi _j=\phi _j
\circ \pi |_{PU_j}$ on $\pi ^{-1}(U_j)$ for each $j$,
where $pr_1: U_j\times H\to U_j$ and $pr _2: U_j\times H\to H$
are projections, $P\phi _j$ is bijective, $P_x\phi _j=P\phi _j|_{H_x}:
H_x\to \{ \phi _j(x) \} \times H$ is a Banach space isomorphism,
$P\phi _l\circ (P\phi _j|_{P(U_l\cap U_j)})^{-1}:
\phi _j(U_j\cap U_l)\times H\to \phi _l(U_j\cap U_l)\times H$
is a $C^n$-diffeomorphism, $l\ge n$. Two atlases are called
equivalent, if their union is an atlas. $(P,M,\pi )$ is called
a vector bundle over $M$ with fibre on $H$. $P$ is called
the total space of $\pi $ and $M$ the base space of $\pi $.
\par Let $(P_1,M_1,\pi _1)$ and $(P_2,M_2,\pi _2)$ be two vector
bundles with spaces $H_1$ and $H_2$ for the fibres of $\pi _1$
and $\pi _2$ respectively. Suppose there are two $C^n$-mappings
$F: M_1\to M_2$ and $PF: P_1\to P_2$ such that $\pi _2\circ PF=
F\circ \pi $ on $P_1$ and the restriction $P_xF:=PF|_{H_{1,x}}:
H_{1,x}\to H_{2,F(x)}$ is a $\bf K$-linear mapping. Then
$(F,PF)$ is called a morphism from $(P_1,M_1,\pi _1)$ to
$(P_2,M_2,\pi _2)$.
\par {\bf 3.1.3.} A $C^m$-vector field on $M$ is a $C^m$-mapping
$\Psi : M\to TM$ such that $\tau _M\circ \Psi =id$.
If $F: M\to N$ is a $C^m$-morphism and $\Psi : M\to TN$
is such that $\tau _N\circ \Psi =F$, then $\Psi $ is called
a vector field along $F$.
\par Suppose that $\bf K$ is spherically complete, then a topologically
adjoint space $H^*$ of $\bf K$-linear functionals on a
Banach space $H$ over $\bf K$ separates points of $H$,
$H^*\ne \emptyset $ (see Lemma 4.3.5 \cite{roo}).
The bundle of $r$-fold contravariant and $s$-fold covariant
tensors over $M$ is defined by $L(\tau ^*,...,\tau ^*,
\tau ,...,\tau ;\rho ): L(T^*M,...,T^*M,TM,...,TM;{\bf K}M)\to M$
or shortly $\tau ^r_s: T^r_sM\to M$, where $\tau ^*$ and
$T^*M$ are repeated $r$ times, $\tau $ and $TM$ are repeated
$s$ times, $\rho : {\bf K}M=M\times {\bf K}\to M$ is the trivial
bundle over $M$. Here $L(\alpha _1,...,\alpha _r;\beta ):
L(A_1,...,A_r;B)\to M$ denotes a vector bundle over $M$,
where $(A_j,M,\alpha _j)$ and $(B,M,\beta )$ are vector bundles,
$\alpha _k^{-1}(x)=H_{k,x}$, $k=1,...,r$, $\beta _H^{-1}(x)=Y_x$,
$L(A_1,...,A_r;B):=\bigcup_{x\in M}L(H_{1,x},...,H_{r,x};Y_x)$,
$H_{k,x}$ and $Y_x$ are isomorphic to Banach spaces $H_k$ and $Y$
respectively over $\bf K$. For each chart $(U_j,\phi _j)$ of $M$
the bundle chart $(U_j,\phi _j,L(A_1,...,A_r;B)\phi _j)$ is given
by $L(A_1\phi _j(x),...,A_r\phi _j(x);B\phi _j(x)):
L(H_{1,x},...,H_{r,x};Y_x)$ $\to L(H_1,...,H_r;Y)$ such that
for $\Psi _x\in L(H_{1,x},...,H_{r,x};Y_x)$ its image
is $L(A_1\phi _j(x),...;B\phi _j(x))\Psi _x=B\phi _j(x)\circ
\Psi _x\circ (A_1\phi _j(x)^{-1}\times ...\times A_r\phi _j(x)^{-1})$,
$A_k\phi _j(x)^{-1}: H_{k,x}\to H_k$ is the $\bf K$-linear isomorphism
of Banach spaces, $L(H_1,...,H_r;Y)$ is the Banach space
of all continuous mappings $f: H_1\times ...\times H_r\to Y$
such that $f$ is $\bf K$-linear by each variable $z_k\in H_k$,
$k=1,...,r$.
\par If $\Psi : M\to \kappa TM$ is a $C^m$-mapping such that
$\kappa \circ \Psi =id$, then $\Psi $ is called a tensor field
(of type $\kappa $), where $(\kappa TM,M,\kappa (\tau ))$
is a tensor bundle over $M$. If $(P,N,\pi )$ is a vector bundle
and $F: M\to N$ is a morphism, then a morphism $\theta :
M\to P$ with $\pi \circ \theta =F$ is called a section along $F$.
\par {\bf 3.1.4.} Let $M$ be a $C^n$-manifold on a Banach space
$X$ over a spherically complete non-Archimedean field $\bf K$
and ${\cal B}_nM$ denotes the set of all $C^n$-vector fields on $M$,
where $n\ge 2$.
Let $\Gamma =\mbox{ }_j\Gamma: \phi _j(U_j)\ni y_j\mapsto
\Gamma (y_j)\in L(X,X;X)$ be a $C^{n-2}$-mapping such that
$$(i)\quad \phi _{l,j}'.\mbox{ }_j\Gamma (y_j)=
\phi _{l,j}"+\mbox{ }_l\Gamma (y_l)
\circ (\phi _{l,j}'\times \phi _{l,j}')$$
for each two charts
with $U_j\cap U_l\ne \emptyset $. This $\{ \mbox{ }_j\Gamma \} $
is called the family of Christoffel symbols $\mbox{ }_j\Gamma $ on $M$.
\par A covariant derivation ${\cal B}_{n-1}M^2\ni (\Psi ,\Phi )\mapsto
\nabla _{\Psi }\Phi \in {\cal B}_{n-2}M$ is given by
$$(ii)\quad \nabla _{\Psi }\Phi (y_j)=
\Phi '(y_j).\Psi (y_j)+\Gamma (y_j)(\Psi (y_j),\Phi (y_j)),$$
where $\Psi (y_j)$ and $\Phi (y_j)$ are principal parts
of $\Psi $ and $\Phi $ on $(U_j,\phi _j)$. If $M$ with $At (M)$
is supplied with $\Gamma $, then $M$ possesses a covariant derivation.
\par {\bf 3.1.5.} For a $C^n$-vector bundle $(P,M,\pi )$ on $X\times H$
with $n\ge 2$ define a $\bf K$-(linear) connection as a bundle morphism
$K: TP\to P$ such that $\pi \circ K=\pi \circ \tau _P$.
This mapping $K$ in its local representation
$\mbox{ }_jK=P\phi _j\circ K\circ TP\phi _j^{-1}$
for bundle charts $(U_j,\phi _j,P\phi _j)$ of $(P,M,\pi )$
and $(TU_j,T\phi _j,TP\phi _j)$ of $(TP,P,\tau _P)$
is given by $\{ U_j,\Xi \} \times (X\times H)\ni
(x,\Psi ,\Phi ,z)\mapsto (x,z+\mbox{ }_j\Gamma (x)(\Phi ,\Psi ))
\in \{ x \} \times H$. The Christoffel symbol
$\mbox{ }_j\Gamma (x): U_j\to L(X,H;H)$ is of class of smoothness $C^{n-2}$.
For it the horizontal space $T_{\Psi h}$ is defined as the kernel
of $K|_{T_{\Psi }P}: T_{\Psi }P \to H_q$, $q=\pi (\Psi )$.
\par For a section $\Psi : M\to P$ in $(P,M,\pi )$ define
the covariant derivation of $\Psi $ in the direction
$\Phi \in T_xM$ by
$$(i)\quad \nabla _{\Phi }\Psi (x)=K\circ T_x\Psi .\Phi .$$
\par {\bf 3.2.} Let $X$ be either a finite dimensional over
a local field $\bf K$ space or of countable type such that a sequence
of subspaces $S_n$ be given with $S_n\subset S_{n+1}$
and $S_n\ne S_{n+1}$ for each $n\in \bf N$,
$cl (\bigcup_nS_n)=X$, a dimension $dim_{\bf K}S_n=:m(n)$
of $S_n$ over $\bf K$ is finite.
Let $U$ be a clopen bounded subset in $S_n$.
Consider an antiderivation operator $P(l,s)$ on
the Banach space $C((t,s-1),U\to {\bf K})$ of functions
$f: U\to \bf K$ with definite partial difference quotients
having continuous extensions (see \S I.2 \cite{lubp2}) and denote
$P(l,s)$ on $U$ by $P_U(l,s)$, where $t\in [0,\infty )$,
$1\le s\in \bf Z$, $l=[t]+1$, $[t]$ is an integer part of $t$.
In particular,
$C((t,0),U\to {\bf K})$ is denoted here by $C^t(U,{\bf K})$.
\par {\bf 3.3. Definition and Note.} Let now $U$ be a clopen
bounded subset
in $X$ with $dim_{\bf K}X=\infty $. For each $f\in C_0((t,s-1),
U\to {\bf K})$ there exists a sequence of cylindrical functions
$f_n$ such that $f=\sum_nf_n$ and $\lim_n \| {\hat f}_n \|_{C((t,s-1),
U_n\to {\bf K})}=0$, where $f_n$ is a cylindrical function on $U$
such that $f_n(x)={\hat f}_n(\pi _nx)$, ${\hat f}_n$
is a function on $U_n:=S_n\cap U$, $\pi _n: X\to S_n$
is a projection on $S_n$. For each $t<\infty $ there exists
$U$ of sufficiently small diameter $\delta $ such that
$\| P_{U_n}(l,s) \| \le 1$ for each $n$, since it is sufficient
to take $\delta ^{|j|+n}/|(j+{\bar u})!| \le 1$ for each $j$
with $|j'|=0,...,l-1$, $j=j'+s'{\bar u}$, $s'\in \{ 0,1,...,s-1 \} $
(see Definition I.2.11 \cite{lubp2}).
For $U$ of $diam (U)$ satisfying such condition define
$P_U(l,s)f:=\sum_nP_{U_n}(l,s)f_n.$
\par For $U$ as above is defined
the space $\mbox{ }_PC_0((t,s),U\to Y):=P_U(l,s)C_0((t,s-1),U\to Y)$,
where $Y$ is a Banach space over $\bf K$.
\par {\bf 3.4. Lemma.} {\it An image $P_U(t,s)(C((t,s-1),U\to Y))$
denoted by $\mbox{ }_PC((t,s),U\to Y)$ is contained in
$C((t,s),U\to Y)$ and does not coincide with the latter space.
The space $\mbox{ }_PC((t,s),U\to Y)$ can be supplied with
a norm denoted by $\| * \|_{U,(t,s),P}$ relative to which it is
complete and $P_U(l,s): (C((t,s-1),U\to Y),
\| * \| _{C((t,s-1),U\to Y)})
\to (\mbox{ }_PC((t,s),U\to Y), \| * \|_{U,(t,s),P})$ is continuous.}
\par {\bf Proof.} Consider at first $dim_{\bf K}X<\infty $.
If $f\in \mbox{ }_PC((t,s),U\to Y)$,
then $\partial ^{\bar u}(P(t,s)f)(x)=f(x)$ for each $x\in U$
(see Corollary I.2.16 \cite{lubp2}). On the other hand, there are
$g\in C((t,s),U\to Y)$ for which $\partial ^{e_j}g(x)=0$
in the notation of Definitions 2.4.1 and 2.11 I.\cite{lubp2}, for example,
locally constant $g$.
\par Let now $X$ may be infinite dimensional, then from taking
the limit of $f_n$ this statement follows in the general case.
Consider an image $P_U(l,s)(B(C((t,s-1),U\to Y),0,1)))=:V$
of the closed ball in $C((t,s),U\to Y)$ containing $0$
and with the unit radius. Let $f\in \mbox{ }_PC((t,s),U\to Y)$,
then there exists $g\in C((t,s-1),U\to Y)$ such that
$P_U(l,s)g=f$. On the other hand,
$\| g \| _{C((t,s-1),U\to Y)}<\infty $ and there exists
a constant $0\ne c\in \bf K$ such that $cg\in B(C((t,s-1),U
\to Y),0,1))$. Therefore, $cf\in V$, since $P_U(l,s)$
is the $\bf K$-linear operator, that is, $V$ is the absorbing
subset. Since the ball $B(C((t,s-1),U\to Y),0,1)$
is $\bf K$-convex, then $V$ is $\bf K$-convex. Evidently, $0\in V$. 
\par Consider a weak topology on $C((t,s),U\to Y)$, then
it induces a weak topology on its $\bf K$-linear subspace
$\mbox{ }_PC((t,s),U\to Y)$. In particular,
each evaluation functional $h_x(f):=f(x)$ is $\bf K$-linear and
continuous on the latter space, where $x\in U$.
In view of Theorem I.2.15 \cite{lubp2} $P_U(l,s)$ is continuous
from $C((t,s-1),U\to Y) \to C((t,s),U\to Y)$.
Therefore, $V$ is bounded relative to the weak topology,
since $U$ is compact and $V$ is bounded relative to a weaker topology
generated by evaluation functionals. 
Let $\eta $ be a Minkowski functional on $\mbox{ }_PC((t,s),U\to Y)$
generated by $V$. It generates a norm in $\mbox{ }_PC((t,s),U\to Y)$
relative to which it is complete. Since $V$ is the unit ball
relative to this norm and $P_U(l,s)^{-1}(V)$ is the unit ball
in $C((t,s-1),U\to Y)$, then $P_U(l,s)$ is continuous
relative to this topology.
\par {\bf 3.5. Note.} In view of Lemma 3.4 Definitions 3.1.1-3.1.5
can be spread on $C_0((t,s))$ and $\mbox{ }_PC_0((t,s))$-manifolds, that is,
$(\phi _{l,j}-id) \in C_0((t,s),W_{l,j}\to X)$ and
$(\phi _{l,j}-id) \in \mbox{ }_PC_0((t,s),W_{l,j}\to X)$
respectively for each charts $U_l$ and $U_j$
with $U_l\cap U_j\ne \emptyset $, where $\phi _j(U_j)$
are bounded clopen subsets in $X$ of sufficiently small
diameter as in \S 3.3 if $X$ is infinite dimensional over $\bf K$.
\par {\bf 3.6. Note.} Consider the space of functions
${\cal F}_{(t,s)}M=C_0((t,s),M\to {\bf K})$, then
\par $\nabla _S(aV+bW)=a\nabla _SV+b\nabla _SW$,
$\nabla _S(fV)=S(f)V+f\nabla _SV$, where $S$, $V$, $W\in {\cal B}_{(t,s)}M$,
${\cal B}_{(t,s)}M$ denotes the set of all $C_0((t,s))$-vector fields
on $M$. Considering the foliation of $M$ and taking the limit we get
for a given chart $(U_j,\phi _j)$:
\par $\nabla _SV(\phi _j)=\sum_k \{ \sum_iS^i(\phi _j)
(\partial V^k/\partial \phi ^i_j)(\phi _j)$ \\
$+\sum_{i,l}S^i(\phi _j)V^l(\phi _j)\Gamma ^k_{i,l}(\phi _j) \} e_k$, \\
where $(\phi _j,e_i)$ are basic vector fields on $\phi _j(U_j)$,
$S(\phi _j)=\sum_iS^i(\phi _j)e_i$, $\Gamma (\phi _j)=
\sum_{i,l,k} \Gamma ^k_{i,l}(\phi _j)e^i\otimes e^j\otimes e_k$,
$e^i(e_j)=\delta ^i_j$ for each $i$ and $j\in \alpha $.
Therefore, there exists a torsion tensor \\
$T(S,V)=\nabla _SV-\nabla _VS-[S,V]$ and a curvature tensor \\
$R(S,V)W=\nabla _S\nabla _VW-\nabla _V\nabla _SW-\nabla _{[S,V]}W$ \\
for each $S$, $V$ and $W\in {\cal B}_{(t,s)}M$ such that 
$T(S,V)=-T(V,S)$, $R(S,V)W=-R(V,S)W$ and \\
$T(\phi _j)(S,V)=\Gamma (\phi _j)(S,V)-\Gamma (\phi _j)(V,S)\in L(X,X;X)$,\\
$R(\phi _j)(S,V)W=D\Gamma (\phi _j).S(V,W)-D\Gamma (\phi _j).V(S,W)$ \\
$+\Gamma (\phi _j)(S,\Gamma (\phi _j)(V,W))-\Gamma (\phi _j)(V,
\Gamma (\phi _j)(S,W))\in L(X,X,X;X)$ analogously to Lemma 1.5.3
\cite{kling}.
\par {\bf 3.7. Theorem.} {\it Let $M$ be a $\mbox{ }_PC_0((t,s))$-manifold
with $s\ge 2$, then there exists a clopen neighbourhood ${\tilde T}M$
of $M$ in $TM$ and an exponential $C_0((t,s))$-mapping
$\exp : {\tilde T}M\to M$ of ${\tilde T}M$ on $M$.}
\par {\bf Proof.} Let $M$ be embedded into $TM$ as the zero
section of the bundle $\tau _M$. Consider the non-Archimedean
geodesic equation $\nabla _{\dot c}{\dot c}=0$
with initial conditions $c(0)=x_0$, ${\dot c}(0)=y_0$,
$x_0\in M$, $y_0\in T_{x_0}M$, where
$c(b)$ is a $\mbox{ }_PC_0((t,s))$-curve on $M$, $c: B({\bf K},0,1)\to M$.
For a chart $(U_j,\phi _j)$ containing a point $x$ of $M$
let $\phi _j\circ c(b):=\psi _j(b)$, thus
$$(i)\quad {\psi _j}"(b)+
\Gamma (\psi _j(b))({\dot \psi }_j(b),{\dot \psi }_j(b))=0.$$
Since $\psi _j\in \mbox{ }_PC_0((t,s))$, then there exists
$f\in C((t,s-2),B\to X)$ such that $\psi _j=
P_B(l,s)P_B(l,s-1)f$, where $B:=B({\bf K},0,1).$
Therefore, ${\dot \psi }_j=P_B(l,s-1)f$ and ${\psi _j}"=f$,
consequently, $f$ satisfies the equation
$$(ii)\quad f(b)+\Gamma (P^2f|_b)(P^1f|_b,P^1f|_b)=0,$$
where $P^2:=P_B(l,s)P_B(l,s-1)$ and $P^1:=P_B(l,s-1)$.
Consider a marked point $b_0\in B$. At first there exists
$r>0$ such that Equation $(ii)$ and hence $(i)$ has a unique
solution in $B({\bf K},b_0,r)$. For this consider the iterational
equation:
$$(iii)\quad f_{m+1}(b)+\Gamma (P^2f_m|_b)(P^1f_m|_b,P^1f_m|_b)=0,$$
where $f_m$ is a sequence of functions.
From $\Gamma \in \mbox{ }_PC_0((t,s-2))$, since $M$ is the
$\mbox{ }_PC_0((t,s))$-manifold, it follows, that $f_{m+1}
\in \mbox{ }_PC((t,s-2))$ for each $f_m\in \mbox{ }_PC((t,s-2))$.
Then $f_{m+1}(t)-f_m(t)=-\Gamma (P^2f_m|_t)(P^1f_m|_t,P^1f_m|_t)+$
$\Gamma (P^2f_m|_t)(P^1f_{m-1}|_t,P^1f_{m-1}|_t)-$
$\Gamma (P^2f_m|_t)(P^1f_{m-1}|_t,P^1f_{m-1}|_t)$ \\
$+\Gamma (P^2f_{m-1}|_t)(P^1f_{m-1}|_t,P^1f_{m-1}|_t)$.
In view of the ultrametric inequality,
bilinearity of $\Gamma (x)(a,b)$ by $a$, $b$
and continuity by $x$, and continuity of $P^1$ and $P^2$
for each $x_0\in M$ and each $t_0\in B({\bf K},0,1)$
there exists $r>0$ and $\epsilon >0$ such that \\
$(iv)$ $\| f_{m+1}-f_m \| \le C \epsilon ^2 \| \Gamma \|
\| f_m-f_{m-1} \| $ for each $t\in B({\bf K},t_0,r)$ and
each $\| y_0 \| <\epsilon $, where $C>0$ is a constant
related with $P^1$ and $P^2$. There exists $0<r<\infty $
such that $\| P^1 \| \le 1$ and $\| P^2 \| \le 1$ and
$P^2f\in G_{j,k}\subset U_j$ for each $f\in G_{j,k}$, since $t$
and $s$ are finite (see above), where $G_{j,k}$ is a clopen
subset in $U_j$, $\| \Gamma \| $ is a norm of
$\Gamma $ on $G_{j,k}\times X^2$ as a bilinear operator on $X$ for
each $x\in G_{j,k}$. In view of continuity of $\Gamma $ and boundedness
of $\phi _j(U_j)$ for each $j$ it is possible
to choose a locally finite covering $G_{j,k}$
subordinated to $U_j$ such that
$\| \Gamma \| $ is finite on $G_{j,k}$, $k \in \bf N$.
Therefore, choosing $C\epsilon ^2 \| \Gamma \| <1$ we get a convergent
sequence on $B({\bf K},t_0,r)\times G_{j,k}\times B(X,0,\delta )$
and due to the fixed point theorem there exists a unique solution
in $B({\bf K},t_0,r)$. In view of compactness of
$B({\bf K},0,1)$ there exists a solution on it.
Let $f$ and $g$ be two functions providing solutions $\psi ^f=P^2f$ and
$\psi ^g=P^2g$ of the problem on $B({\bf K},0,1)$, then
$P^2f(t_l)=P^2g(t_l)$, $P^1f(t_l)=P^1g(t_l)$ for a finite
number of points $t_0=0,$ $t_1,...,t_k\in B({\bf K},0,1)$
such that on each $B({\bf K},t_j,r_j)$ a solution is unique
for a given initial conditions, $0<r_j\le 1$ for each $j$
and $\bigcup_jB({\bf K},t_j,r_j)=B({\bf K},0,1)$. This imply that \\
$(v)$ $\Gamma (P^2f|_t)(P^1f|_t,P^1f|_t)-
\Gamma (P^2(f+c_{1,l})|_t)(P^1(f+c_{2,l})|_t,P^1(f+c_{2,l})|_t)=c_{1,l}$ \\
for each $l$ and each $t\in B({\bf K},t_l,r_l)$.
On the other hand, $P^1c$ and $P^2c$ are not locally constant for
a constant $c\ne 0$, $\Gamma (\phi _j)(a,b)$ is bilinear
by $(a,b)\in X^2$ and satisfies Equation $3.1.4.(i)$,
hence Equation $(v)$ may be satisfied
only for $c_{1,l}=c_{2,l}=0$ for each $l$, consequently,
a solution is unique.
\par Since $f\in \mbox{ }_PC_0((t,s-2))$, then $\psi _j\in \mbox{ }_P
C_0((t,s))$ for each $j$.  Moreover, $c_{aS}(t)=c_S(at)$
for each $a\in B({\bf K},0,1)$ such that $|aS(\phi _j(q))|<\epsilon ,$
since $dc_S(at)/dt=a(dc_S(z)/dz)|_{z=at}$.
In view of continuity of $P^2$ and $P^1$ and $\Gamma $ operators,
for each $x_0\in M$ there exists a chart $(U_j,\phi _j)$
and clopen neighbourhoods $V_1$ and $V_2$,
$\phi _j(x_0)\in V_1\subset V_2\subset
\phi _j(U_j)$ and $\delta >0$ such that from $S\in TM$ with
$\tau _MS=q\in \phi _j^{-1}(V_1)$ and $|S(\phi _j(q))|<\delta $
it follows, that the geodesic $c_S$ with $c_S(0)=S$ is defined
for each $t\in B({\bf K},0,1)$ and $c_S(t)\in \phi _j^{-1}(V_2)$.
Due to paracompactness of $TM$ and $M$ this covering
can be chosen locally finite \cite{eng}.
\par This means that there exists a clopen neighbourhood
${\tilde T}M$ of $M$ in $TM$ such that a geodesic $c_S(t)$
is defined for each $S\in {\tilde T}M$ and each $t\in B({\bf K},0,1)$.
Therefore, define the exponential mapping $\exp : {\tilde T}M\to M$
by $S\mapsto c_S(1)$, denote by $\exp _x:=\exp |_{{\tilde T}M\cap T_xM}$
a restriction to a fibre. Then $\exp $ has a local representation
$(x_0,y_0)\in V_1\times B(X,0,\delta )\mapsto \psi _j(1;x_0,y_0)\in
V_2\subset \phi _j(U_j).$ From Equations $(i,ii)$ it follows that
$\exp $ is of $C_0((t,s))$-class of smoothness from ${\tilde T}M$
onto $M$.
\par {\bf 3.8. Corollary.} {\it If $M$ is a $\mbox{ }_PC_0((t,s))\cap
C^{\infty }$-manifold with $s\ge 2$, then $\exp \in C^{\infty }
({\tilde T}M,M).$}
\par {\bf 3.9. Note.} If $M$ is an analytic manifold, then
$\exp : {\tilde T}M\to M$ is a locally analytic mapping.
Theorem 3.7 gives an exponential manifold mapping for wider
class of manifolds, than treated by the rigid geometry.
\par {\bf 3.10. Note and Definitions.}
Let $M$ be a $C^{\infty }$-manifold and let $\tau _M: TM\to M$
be the tangent bundle, $\theta : M\times H\to M$ be a trivial
bundle over $M$ with a Banach space fibre $H$ over $\bf K$.
There exists the bundle $L_{1,r}(\theta ,\tau _M)$ over $M$ with
the fibre $L_{1,r}(H,X)$, where $r\ge 1$ and spaces
$L_{n,r}(H,X)$ were defined in \cite{lufpmsp} and \S 2 in
I.\cite{lustpr} (but with the notation $L_{r,n}$ there).
\par Let $M$ be a $C^{\infty }$-manifold with functions
$\phi _{l,j}$ satisfying Conditions I.4.8.(i) \cite{lustpr}.
Suppose that $w$ is a stochastic process with values
in $H$ and $\xi $ is a stochastic process with values in $X$
such that $\lambda \{ \omega : w(t,\omega )\in C^0\setminus
C^1 \} =0$, where $H$ and $X$ are Banach spaces over a
local field $\bf K$. \\
Let $a\in L^q(\Omega ,{\sf F},\lambda ;C^0(B_R,L^q(\Omega ,
{\sf F},\lambda ;C^0(B_R,X))))$ and \\
$E\in L^r(\Omega ,{\sf F},\lambda ;C^0(B_R,L(
L^q(\Omega ,{\sf F},\lambda ;C^0(B_R,H),
L^q(\Omega ,{\sf F},\lambda ;C^0(B_R,X)))))$,  \\
$(i) \quad \xi (t,\omega )=
\xi _0(\omega )+({\hat P}_ua)(u,\omega ,\xi )|_{u=t}+
({\hat P}_{w(u,\omega )}E)(u,\omega ,\xi )|_{u=t}$, \\
where $1\le r, s, q\le \infty $, $1/r+1/s=1/q$,
$w\in L^s(\Omega ,{\sf F},\lambda ;C^0_0(B_R,H))$,
$\xi \in L^q(\Omega ,{\sf F},\lambda ;C^0(B_R,X))$.
Since $H$ and $X$ are isomorphic with $c_0(\alpha _H,{\bf K})$
and $c_0(\alpha _X,{\bf K})$, then $L_{n,r}(X,H)$ has the embedding
into $L_{n,r}(H,H)$ for $\alpha _X\subset \alpha _H$
and $L_{n,r}(H,H)$ has an embedding into $L_{n,r}(X,X)$
for $\alpha _H\subset \alpha _X$. Inclusions
$Range (E)\subset X$, $Range (w)\subset H$ and
$Range (\xi )\subset X$ reduce this case to Theorem II.3.4 \cite{lustpr}.
In view of Lemma $I.2.3$ and Formula $I.4.8.(ii)$ \cite{lustpr}
$$(ii)\quad d\phi (\xi (t,\omega ))=J(\phi ,a,E)adt+J(\phi ,a,E)Edw,
\mbox{ where},$$
$$(iii)\quad J(\phi ,a,E):=\sum_{m=0}^{\infty }[m!]^{-1}
\sum_{l=0}^m{m\choose l}{\hat P}_{u^l,w^{m-l}}\phi ^{(m+1)}
\circ (a^{\otimes l}\otimes E^{\otimes (m-l)})|_{u=t}.$$
\par For stochastic processes of type II.3.5 \cite{lustpr}
it is necessary to consider the following generalization of Theorem
I.4.8 \cite{lustpr}.
\par {\bf 3.11. Note.} Consider $a\in 
L^{\infty }(\Omega ,{\sf F},\lambda ;C^0(B_R,L^q(
\Omega ,{\sf F},\lambda ;C^0(B_R,X))))$
and $E\in L^{\infty }(\Omega ,{\sf F},\lambda ;C^0(B_R,L(L^q(
\Omega ,{\sf F},\lambda ;C^0(B_R,H),L^q(
\Omega ,{\sf F},\lambda ;C^0(B_R,X))))),$
$a=a(t,\omega ,\xi )$, $E=E(t,\omega ,\xi )$, 
$t\in B_R,$ $\omega \in \Omega ,$
$\xi \in L^q(\Omega ,{\sf F},\lambda ;C^0(B_R,X))$
and $\xi _0 \in L^q(\Omega ,{\sf F},\lambda ;X),$
$w\in L^{\infty }(\Omega ,{\sf F},\lambda ;C^0_0(B_R,H)),$
$1\le q\le \infty $, 
where $a$ and $E$ satisfy the local Lipschitz condition
(see II.3.4.(LLC)\cite{lustpr}).
Suppose $\xi $ is a stochastic process of the type
\par $(i)$ $\xi (t,\omega )=\xi _0(\omega )+$
$\sum_{m+b=1}^{\infty }\sum_{l=0}^m(
{\hat P}_{u^{b+m-l},w(u,\omega )^l}[a_{m-l+b,l}
(u,\xi (u,\omega ))\circ (I^{\otimes b}\otimes 
a^{\otimes (m-l)}\otimes E^{\otimes l})])
|_{u=t}$  \\
such that
$a_{m-l,l} \in C^0(B_{R_1}\times B(L^q(\Omega ,{\sf F},\lambda ;C^0(B_R,X)),
0,R_2),L_m(X^{\otimes m};X))$ (continuous and bounded on its domain)
for each $n, l,$ $0<R_2<\infty $ and
\par $(ii)$ $\lim_{n\to \infty } \sup_{0\le l\le n}\|a_{n-l,l}
\|_{C^0(B_{R_1}\times B(L^q(\Omega ,{\sf F},\lambda ;C^0(B_R,X)),
0,R_2),L_n(X^{\otimes n},X))}=0$ 
for each $0<R_1\le R$ when $0<R<\infty $, or each $0<R_1<R$
when $R=\infty $, for each $0<R_2<\infty .$
\par Moreover, suppose that a function $f$ satisfies the conditions:
$$(iii)\quad f(u,x)\in C^{\infty }(T\times H,X)$$ 
and 
$$(iv)\mbox{ } \lim_{n\to \infty }
\max_{0\le l\le n}\|(\bar \Phi ^nf)(t,x;h_1,...,h_n;$$
$$\zeta _1,...,\zeta _n)
\|_{C^0(T\times B({\bf K},0,r)^l\times B(H,0,1)^{n-l}
\times B({\bf K},0,R_1)^{n-l},X)}=0$$
for each $0<R_1<\infty ,$ where $h_j=e_1$ and
$\zeta _j\in B({\bf K},0,r)$ for variables corresponding to 
$t\in T=B({\bf K},t_0,r)$
and $h_j\in B(H,0,1)$, $\zeta _j\in B({\bf K},0,R_1)$
for variables corresponding to $x\in H$.
\par Analogously $a$, $E$, $a_{l,m}$ for $\xi $ with values
in $M$ are considered substituting $C^0(B_R,H)$ on $C^0(B_R,M)$.
\par {\bf 3.12. Theorem.} {\it If Conditions $3.11.(ii)$
are satisfied, then
$(i)$ has the unique solution in $B_R$.
If in addition Conditions $3.11.(iii,iv)$ are satisfied, then
$$(i)\quad f(t,\xi (t,\omega ))=f(t_0,\xi _0)+
\sum_{m+b\ge 1, 0\le m\in {\bf Z}, 0\le b\in {\bf Z}}
((m+b)!)^{-1}\sum_{l_1,...,l_m}{{m+b}\choose m}$$
$$({\hat P}_{u^{b+m-l},w(u,\omega )^l}[(\partial ^{(m+b)}f/
\partial u^b\partial x^m)
(u,\xi (u,\omega ))\circ (a_{l_1,n_1}\otimes ...\otimes
a_{l_m,n_m})\circ (I^{\otimes b}\otimes
a^{\otimes (m-l)}\otimes E^{\otimes l})])|_{u=t},$$
where $l_1+...+l_m=m+b-l$, $n_1+...+n_m=l$, $l_1,...,l_m,
n_1,...,n_m$ are nonnegative integers.}
\par {\bf Proof.} The first part of the theorem follows
from $II.3.5$ \cite{lustpr} and embeddings of \S 3.10.
Since $\sigma _n\circ \sigma _m(t)=
\sigma _n\circ \sigma _{m+j}(t)$ for each $n\ge m$,
$j>0$ and $\sigma _0(t)=t_0$, then
from Formula $I.2.1.(4)$ it follows, that \\
${\hat P}_{u^{l+b},w^m}a_{l+b,m}\circ (I^{\otimes b}\otimes a^{\otimes l}
\otimes E^{\otimes m})|^{u=t_{n+1}}_{u=t_n}=a_{l+b,m}(t_n)\circ
((t_{n+1}-t_n)^{\otimes b}\otimes (a(t_n)(t_{n+1}-t_n))^{\otimes l}
\otimes (E(t_n)(w(t_{n+1}-w(t_n)))^{\otimes m}),$ \\
where other arguments are omitted for shortening the notation.
Therefore, the second part of this theorem follows from Formulas
$I.4.8.(iii)$ \cite{lustpr} and $3.11.(i)$.
\par {\bf 3.13. Note.} Let Conditions $3.11.(i-iv)$ be satisfied
and $\phi =f$ be independent from $t$.
Then due to Lemma $I.2.3$ \cite{lustpr} and Theorem $3.12$
above Formula $3.10.(ii)$ is valid with new operator $J$:
$$(i)\quad J(\phi ,a,E):=\sum_{m=0}^{\infty }[m!]^{-1}
\sum_{l_1,...,l_m}{\hat P}_{u^l,w^{m-l}}\phi ^{(m+1)}
\circ (a_{l_1,n_1}...a_{l_m,n_m})\circ
(a^{\otimes l}\otimes E^{\otimes (m-l)}),$$
where $l_1+...+l_m=l$, $n_1+...+n_m=m-l$.
\par {\bf 3.14. Definition.} Let $(\Pi ,M,\pi )$
be a bundle on a manifold $M$ with fibres $X\oplus L(H,X)$
for each $x\in M$ and with transition functions $J(\phi ,a,E):
(a,E)\mapsto (J(\phi ,a,E)a,J(\phi ,a,E)E)$,
where $\phi =\phi _{j,l}$ for each pair of charts
$(U_j,\phi _j)$ and $(U_l,\phi _l)$ with $U_j\cap U_l\ne \emptyset $,
$a\in X$, $E\in L(H,X)$, $J(\phi ,a,E)$ is given either
by Equation $3.10.(iii)$ or by $3.13.(i)$.
\par {\bf 3.15. Definition and Note.} Let $t\in T\subset \bf K$,
where $\bf K$ is a local field, $T$ is clopen in $\bf K$.
Let also $(U_j,\phi _j)$ be a chart of a manifold $M$ on a Banach space
$X$ over $\bf K$, $x\in U_j\subset M$, $(a,E)\in \pi ^{-1}(x)$
(see \S 3.14). By ${\cal G}_x(a,E)$ is denoted a collection
of $M$-valued stochastic processes $\xi $ such that $\xi \in U_j$
with probability $1$, where $\phi _j\circ \xi $ is a solution of
Equation either $3.10.(i)$ or $3.11.(i)$ for each $j$.
Then ${\cal G}_x(a,E)$
is called the germ of the diffusion process at the point
$x$ defined by a pair $(a,E)$. It is in addition with a given family of
sections $a_{l,m}$ of bundles $(\Pi _{l+m},M,\pi _{l+m})$
with fibres $L_{m+l}(X^{\otimes {m+l}};X)$
such that $a_{l,m,x}\in \pi _{l+m}^{-1}(x)$ in the case $3.11$.
Therefore, $3.10$ is the particular case of $3.11$.
\par A section $\cal U$ of the vector bundle $(\Pi ,M,\pi )$
is the non-Archimedean analog of It$\hat o$'s field over $M$.
\par {\bf 3.16. Theorem.} {\it Let $\phi $ and $\psi $ be two functions
satisfying conditions either of \S 3.10 or \S 3.11 such that
$Dom (\phi )\supset Range (\psi )$. Then
\par $(i)$ $J_{\psi (x)}(\phi ,a,E)\circ J_x(\psi ,a,E)=
J_x(\phi \circ \psi ,a,E)$,
\par $(ii)$ $J_x(id,a,E)=id$.}
\par {\bf Proof.} Since $a_{l,m,x}\in L_{l+m}(X^{\otimes {l+m}};X)$,
then $J_x(\phi ,a,E)a_{l,m,x}\circ (a^{\otimes l}\otimes E^{\otimes m})=
a_{l,m,x}\circ ((J_x(\phi ,a,E)a)^{\otimes l}\otimes
(J_x(\phi ,a,E)E)^{\otimes m})$ for each $0\le l, m\in \bf Z$
and $x\in M$, where $a=a_x$, $E=E_x$, $(a_x,E_x)\in \pi ^{-1}(x)$.
Each derivative $\phi ^{(m)}$ and $\psi ^{(m)}$ is a $m$-polylinear
operator on $X$. Therefore,
$(\phi \circ \psi )^{(m)}(x)=\sum_{l_1+...+l_b\ge m, 1\le b\le m}
R_b\circ (Q_{l_1}\otimes ...\otimes Q_{l_b})$, where
$R_b$ and $Q_l$ are the $b$-linear and $l$-linear operators
corresponding up to constant multipliers
to $\phi ^{(b)}(z)|_{z=\psi (x)}$ and $\psi ^{(l)}(x)$.
Then $\sum_kQ_{l_j}(\Delta _k\xi _1,...,\Delta _k\xi _{l_j})=
{\hat P}_{u^{l_{j,1}},w^{l_{j,2}}}Q_{l_j}(a^{\otimes l_{j,1}}\otimes
E^{\otimes l_{j,2}})$ for nonegative $l_{j,1}$ and $l_{j,2}$
with $l_{j,1}+l_{j,2}=l_j$ and $\xi _i(t,\omega )={\hat P}_ua|_{u=t}$
for $i=1,...,l_{j,1}$, $\xi _i(t,\omega )={\hat P}_wE|_{u=t}$
for $i=l_{j,1}+1,...,l_j$. Moreover, \\
$\sum_kQ_{l_j}(\Delta _k\xi _1,...,\Delta _k\xi _{l_j-1},\xi _{l_j})=
{\hat P}_{u^{l_{j,1}},w^{l_{j,2}}}Q_{l_j}(a^{\otimes l_{j,1}}\otimes
E^{\otimes l_{j,2}})v$ \\
for nonegative $l_{j,1}$ and $l_{j,2}$
$l_{j,1}+l_{j,2}=l_j-1$ and $\xi _i(t,\omega )={\hat P}_ua|_{u=t}$
for $i=1,...,l_{j,1}$, $\xi _i(t,\omega )={\hat P}_wE|_{u=t}$
for $i=l_{j,1}+1,...,l_j-1$, $\xi _{l_j}=v$, where either $v=a$
or $v=E$.
Therefore, $\phi : {\cal G}_x(a,E)\to {\cal G}_{\phi (x)}(Ja,JE)$,
where $J=J(\phi ,a,E)$. In view of Theorems $I.4.8$ \cite{lustpr}
and $3.12$, Formulas $3.10.(iii)$ and $3.13.(i)$
there is satisfied the equality
\par $(iii)$ $J_x(\phi ,a,E)=\phi '(\xi ^0_x)$, \\
where $\xi ^0$ is a stochastic process being the solution either of
Equation $3.10.(i)$ or $3.11.(i)$, $x\in M$,
$\xi ^0_x\in T_xM$. On the other hand,
$(\phi (\psi )(x))'=\phi '(\psi (x)).\psi '(x)$
for each $x\in Dom (\psi )$. Therefore, from $\xi \in Dom (\psi )$
and $Range (\psi )\subset Dom (\phi )$, Formula
$3.16.(i)$ follows. From $id'=I$, where $I$ is a unit operator,
Formula $3.16.(ii)$ follows.
\par {\bf 3.17. Remark and Definition.}  Apart from the classical case
here the bundle associated with the operator $J(\phi ,a,E)$
in general is not quadratic. It may be polynomial only
in a particular case given by theorem $I.4.6$ \cite{lustpr}.
\par Let $f=\exp $, where $\exp :=\exp ^M$ is the exponential
mapping for $M$. Consider ${\cal G}_{(x,0)}(a,E)$
a stochastic processes germ at a point $y=0$ in the tangent space
$T_xM$. Then $\exp ^*_x{\cal G}_{(x,0)}(a,E)={\cal G}_x(J(\exp _x,a,E))(a,E)$
is a stochastic processes germ at $x\in M$.
Therefore, $\phi _j\circ \exp ^*_x{\cal G}_{(x,0)}(a,E)=
{\cal G}_{\phi (x)}(J(\phi _j\circ \exp _x,a,E))(a,E)$
for each chart $(U_j,\phi _j)$ of $M$. The germ $\exp ^*_x{\cal G}_{(x,0)}
(a,E)$ is called a stochastic differential bundle.
\par {\bf 3.18. Corollary.} {\it To a functor $J$ there corresponds
a bundle $(J^M,M,\pi _J)$ and a fibre $J^M_x:=\pi _J^{-1}(x)$
may be identified with the space ${\cal G}_x(J^M)$ of stochastic processes
germ. To a morphism $f: M\to N$ of manifolds there corresponds
a bundle morphism ${\cal G}(f)=f*f^*$, where $f^*\xi :=f(\xi )$.}
\par {\bf Proof.} If $f: M\to N$ is a manifold morphism,
then $\cal U$ is transformed as $(a_x,E_x)\mapsto (a^f_{f(x)},
E^f_{f(x)})$, where $a^f_{f(x)}=J(f,a,E)a_x$
and $E^f_{f(x)}=J(f,a,E)E_x$,
$a^f_{l,m,f(x)}(t,f^*\xi )=a_{l,m,x}(t,\xi )$ for each $x\in M$.
The stochastic process $\xi ^0_x$ satisfies the antiderivational equation
$$(i)\quad \xi ^0_x=\sum_{l,m} {\hat P}_{u^l,w^m}
a_{l,m,x}\circ (a_x^{\otimes l}\otimes E_x^{\otimes m})$$
and its differential has the form:
$$(ii)\quad d\xi ^0_x=\sum_{l,m}l{\hat P}_{u^{l-1},w^m}
a_{l,m,x}\circ (a_x^{\otimes (l-1)}\otimes E_x^{\otimes m})a_xdt$$
$$+\sum_{l,m}m{\hat P}_{u^l,w^{m-1}}
a_{l,m,x}\circ (a_x^{\otimes l}\otimes E_x^{\otimes (m-1)})E_xdw_x.$$
Hence
$$(iii)\quad f(\xi ^0_x(t,\omega ))=\sum_{l,m} {\hat P}_{u^l,w^m}
a^f_{l,m,f(x)}(u,f(\xi ^0_x(u,\omega )))\circ ({a^f_{f(x)}}^{\otimes l}
\otimes {E^f_{f(x)}}^{\otimes m})|_{u=t}.$$
Therefore, $f^*: \exp_x^{M*} (d\xi ^0_x)\mapsto
\exp_{f(x)}^{N*}(f^*d\xi ^0_x),$ where
$$(iv)\quad f^* d\xi ^0_x=\sum_{l,m}l{\hat P}_{u^{l-1},w^m}
a^f_{l,m,f(x)}\circ ({a^f_{f(x)}}^{\otimes (l-1)}\otimes {E^f_{f(x)}}^{
\otimes m})a^f_{f(x)}dt$$
$$+\sum_{l,m}m{\hat P}_{u^l,w^{m-1}}
a^f_{l,m,f(x)}\circ ({a^f_{f(x)}}^{\otimes l}\otimes {E^f_{f(x)}}^{
\otimes (m-1)}){E^f_{f(x)}}df(w_x)$$ for $f$-related mappings
$\exp^M$ and $\exp^N$.
\par {\bf 3.19. Theorem.} {\it Let $\exp : {\tilde T}M\to M$ be
the exponential mapping of a manifold $M$. Then $J(\exp ,a,E):
J^{{\tilde T}M}\to J^M$ is a bundle morphism. If $(U,\phi )$ is a chart
of $M$, then
$$(i)\quad J(\exp ^{\phi }_{\phi (x)},a_x^{\phi },E_x^{\phi })
(a_x^{\phi },E_x^{\phi })=
(Sa_x^{\phi },SE_x^{\phi }),$$ where 
$S:=(d[\phi \circ \exp _x\circ [\phi '(x)]^{-1}](z)/dz)|_{z=\xi ^0_x}.$}
\par {\bf Proof.} The first statement of the theorem follows from
Theorem 3.16 and Corollary 3.18. Consider the mapping
$F(z):=[\phi \circ \exp _x\circ [\phi '(x)]^{-1}](z)$
for a chart $(TU_j,T\phi _j)$ of $TM$, where $\phi =\phi _j$.
The mapping $F$ is the local representation of $\exp $
in terms of coordinate mappings.
Hence $J(\exp ^{\phi }_{\phi (x)},a_x^{\phi },E_x^{\phi })=
[dF(z)/dz]|_{z=\xi ^0_x}$, where $\xi ^0_x$
is a solution of Equation either $3.10(i)$ or $3.11.(i)$ in $T_xM$.
In particular, $F'(0)=id$ and $F"(0).(v,v)=-\Gamma (x)(v,v)$,
but in general $\xi ^0_x$ may be nonzero.
\par {\bf 3.20. Definition.} Let $\cal U$ be a section
of the bundle $(\Pi ,M,\pi )$. Consider the differential
$$(i)\quad d\xi (t,\omega )=
\exp_{\xi (t,\omega )}^* {\cal G}(a_{\xi (t,\omega )},
E_{\xi (t,\omega )})$$
and the corresponding antiderivational equation:
$$(ii)\quad \xi (t,\omega )=\exp_{\xi (t,\omega )}
\{ \sum_{l,b,m} {\hat P}_{u^{l+b},w^m}
a_{b+l,m,\xi (t,\omega )}(u,\xi (u,\omega )))\circ (I^{\otimes b}\otimes
a_{\xi (t,\omega )}^{\otimes l}\otimes E_{\xi (t,\omega )
}^{\otimes m})|_{u=t} \} .$$
Suppose that there exists a neighbourhood $V_x\ni x$ and a
stochastic process belonging to the germ
$\exp_x ({\cal G}(a_x,E_x))={\cal G}_x(J(\exp _x,a_x,E_x))(a_x,E_x)$
such that $P_{s,x} \{ \omega : \xi _x(t,\omega )\in V_x, t\ne s \} =1,$
where $P_{s,x}(W):=P(W: \xi (s,\omega )=x)$, $W\in \sf F$.
If this is satisfied for $\nu _{\xi (s)}$-a.e. $x\in M$,
then it is said, that $\xi (t,\omega )$ possesses a stochastic differential
governed by the field $\cal U$, where $\nu _{\xi (s)}(*):=P\circ \xi ^{-1}
(s,*)$. An $M$-valued $\xi $ satisfying $(ii)$
is called an integral process of the field ${\cal U}(t)$.
\par {\bf 3.21. Definition.} An atlas $At(M)=\{ (U_j,\phi _j): j \} $ 
of a manifold $M$ on a Banach space $X$ over $\bf K$
is called uniform, if its charts satisfy 
the following conditions: \\
$(U1)$ for each $x\in M$ there exist
neighbourhoods $U_x^2\subset U_x^1\subset U_j$
such that for each $y\in U_x^2$ there is the inclusion
$U_x^2\subset U_y^1$; \\
$(U2)$ the image $\phi _j(U_x^2)\subset X$
contains a ball of the fixed positive radius
$\phi _j(U_x^2)\supset B(X,0,r):=\{ y: y\in X, \| y\| \le r \} ;$ \\
$(U3)$ for each pair of intersecting charts
$(U_1,\phi _1)$ and $(U_2,\phi _2)$ transition mappings
$\phi _{l,j}=\phi _l\circ \phi _j^{-1}$ are such that
$\sup_x \| \phi _{l,j}' \| \le  C$ and $\sup_x \| \phi _{l,j}(x) \|
\le  C$, where $C=const >0$ does not depend on $\phi _l$ and $\phi _j$.
\par {\bf 3.22. Remark.} Consider a measurable space
$(M,{\cal L})$, where $\cal L$ is a $\sigma $-algebra on $M$,
define a random mapping $S(t,\tau ;\omega ): M\to M$ for each
$t, \tau \in T$ by $x\mapsto S(t,\tau ;\omega ;x)=S(t,\tau ;\omega )
\circ x$. Suppose that
\par $(1)$ the mapping $x\times \omega \mapsto
S(t,\tau ;\omega ;x)$ is ${\cal L}\times \sf F$-measurable
for each $t, \tau \in T$;
\par $(2)$ the random variable $S(t,\tau ;\omega ;x)$ is
$\sf F$-measurable and does not depend on $\sf F$ for each
$t, \tau ,x$, moreover, all others conditions of \S II.3.9 \cite{lustpr}
let also be satisfied (with the notation $S(t,\tau ;\omega )$
here instead of $T(t,s;\omega )$ there).
\par {\bf 3.23. Proposition.} {\it Let $\xi $ be a stochastic
process given by Equation $3.11.(i)$
and let also $\max (\| a(t,\omega ,x)-a(v,\omega ,x) \| ,
\| E(t,\omega ,x)-E(v,\omega ,x) \| ) \le |t-v|(C_1+C_2 \| x \| ^b)$
for each $t$ and $v\in B({\bf K},t_0,R)$ $\lambda $-almost everywhere
by $\omega \in \Omega $,
where $b$, $C_1$ and $C_2$ are non-negative constants.
Then $\xi $ with the probability $1$ has a $C^0$-modification
and  $q(t)\le \max \{ M \| \xi _0 \| ^s,
|t-t_0|(C_1+C_2q(t)) \} $
for each $t\in B({\bf K},t_0,R)$, where
$q(t):= \sup_{|u-t_0|\le |t-t_0|}M \| \xi (u,\omega )\| ^s $
and ${\bf N}\ni s\ge b\ge 0$.
Moreover, if $\lambda \{ \omega : w(t,\omega )\in C^0\setminus C^1 \} =0$,
then for $\lambda $-a.e. $\omega $ there exists $\xi '$ and
$\lambda \{ \omega : \xi (t,\omega )\in C^0\setminus C^1 \} =0$.}
\par {\bf Proof.} In view of Theorem 3.12 applied to
$f(t,x)=x^s$ we have
$$f(t, \xi (t,\omega ))=f(t_0,\xi _0)+\sum_{k=1}^s\sum_{l_1,...,l_k}(
\hat P_{u^{k-l},w(u,\omega )^l}[({s\choose k}\xi (t,\omega )^{s-k}
(u,\xi (u,\omega ))\circ $$
$$(a_{l_1,n_1}\otimes ...\otimes
a_{l_k,n_k})\circ (a^{\otimes (k-l)}\otimes E^{\otimes l})])
|_{u=t},$$ where $l_1+...+l_k=k-l$, $n_1+...+n_k=l$.
From Conditions of \S 3.11 and in particular $3.11.(ii)$
it follows, that
$M \| \xi (t,\omega ) \| ^s\le 
\max (M \| \xi _0 \| ^s,  |t-t_0| 
d(\hat P^s_*) (C_1+C_2 \sup_{|u-t_0|\le |t-t_0|}
M \| \xi (u,\omega )\| ^s)$, since $|t_j-t_0|
\le |t-t_0|$ for each $j\in \bf N$ and
$M \| \xi (t,\omega ) -\xi (v,\omega )\| ^s\le 
|t-v| (1+C_1+C_2 d(\hat P^s_*) \sup_{|u-t_0|\le \max( |t-t_0|, |v-t_0|)}
M \| \xi (u,\omega )\| ^s)$, since $|t_j-v_j|\le |t-v|+\rho ^j$
for each $j\in \bf N$, where $0<\rho <1$, 
$$d( \hat P^s_* ) :=\sup_{(a\ne 0, E\ne 0, f \ne 0, a_{l_j,n_j}\ne 0,
j=1,...,k)} \max_{s\ge k\ge l\ge 0}
\| (k!)^{-1} {\hat P}_{u^{k-l},w^l}(\partial ^kf
/\partial ^kx ) \circ (a_{l_1,n_1}\otimes ... \otimes a_{l_k,n_k})$$
$$\circ (a^{\otimes (k-l)}\otimes E^{\otimes l}) \| /
(\| a \|_{C^0(B_R,H)}^{k-l} \| E \|_{C^0(B_R,L(H))}^l
\| f \|_{C^s(B_R,H)} \prod_{j=1}^k \| a_{l_j,n_j} \| ),$$ 
hence $d( \hat P^s_* )\le 1$, since
$f\in C^s$ as a function by $x$ and $(\bar \Phi ^sg)(x;h_1,...,h_s;0,...,0)=
D^s_xg(x).(h_1,...,h_s)/s!$ for each $g\in C^s$ and due to 
the definition of $\| g \| _{C^s}$. Considering in particular 
polyhomogeneous $g$ on which $d(\hat P^s_*)$ takes its maximum value 
we get $d(\hat P^s_*)=1$.
From conditions on $w$, $a_{l,k}$, $a$ and $E$ it follows, that
$\xi (t,\omega )$ with the probability $1$ has a $C^0$-modification
(see Theorem $II.3.5$ \cite{lustpr}),
$\xi \in L^q(\Omega ,{\sf F},\lambda ;C^0(B_R,H))$.
\par The last statement of this proposition follows from
Lemma $I.2.3$ \cite{lustpr}.
\par {\bf 3.24. Theorem.} {\it Suppose that $M$ is a manifold either
satisfying conditions of Corollary $3.8$ or is analytic, $At (M)$
is uniform (see \S 3.21). Let $a,$ $E,$ $a_{m,l}$ and $w$
corresponding to a section $\cal U$ satisfy
conditions of \S 3.11 with $\lambda \{ \omega : w(t,\omega )\in C^0\setminus
C^1 \} =0$. Then there exists a unique up to stochastic
equivalence random evolution family $S(t,\tau ;\omega )$ for a solution
$\xi (t,\omega )$ of Equation $3.20.(ii)$.}
\par {\bf Proof.} Consider a solution of the
non-Archimedean stochastic equation:
$$(i)\quad \xi (t,\omega )=\exp _{\xi (t,\omega )}  \{
\sum_{m,b,l} {\hat P}_{u^{m+b},w^l}a_{m+b,l,\xi (t,\omega )}
(u,\xi (u,\omega )) \circ (I^{\otimes b}\otimes
a_{\xi (t,\omega )}^{\otimes m}\otimes E_{\xi (t,\omega )}
^{\otimes l})|_{u=t} \} $$
corresponding to $3.20.(i)$.
On each chart of the uniform atlas $At(M)$ of $M$
fields $\{ a_{m,l}: m,l \} $, $a$, $E$ and $w$ are $\lambda $-a.e. bounded
due to conditons of $\S 3.11$.
For each two charts $(U_j,\phi _j)$ and $(U_l,\phi _l)$ with
$U_j\cap U_l\ne \emptyset $ a transition mapping $\phi :=\phi _{j,l}$
is bounded together with its derivatives, hence $\Gamma $ is bounded
on each $U_j$, since the covering $\{ U_j: j \} $ of $M$ can be chosen
locally finite due to paracompactness of $M$ \cite{eng}.
\par In view of Theorem $II.3.5$ \cite{lustpr},
Corollary 3.18 and Theorem 3.19 Equation $(i)$ has a unique solution
on $M$. Let $(a,E)$ be a section of the bundle $(\Pi ,M,\pi )$ and
$a_{l,m}$ be sections of the bundles $(\Pi _{l+m},M,\pi _{l+m})$
(see \S \S 3.14 and 3.15). Consider a family $\zeta _y(x)$
of functions on $M$ of the class $C^1(M,{\bf K})$ such that
$\zeta _y(x)=0$ if $x\notin U^1_y$, $\zeta _y(x)=1$ if $x\in U^2_y$
of the uniform atlas (see \S 3.21), then $a^y_x:=\zeta _y(x)a_x$,
$E^y_x:=\zeta _y(x)E_x$, $a_{l,m,x}^y:=\zeta _y(x) a_{l,m}$
are local fields. Then there exists the local evolution family
$S_y(t,\tau ;\omega )$ for each local solution (that is,
with local coefficients): 
$$(ii)\quad \xi ^y(t,\omega )=\exp _{\xi ^y(t,\omega )}  \{
\sum_{m,b,l} {\hat P}_{u^{m+b},w^l}a^y_{m+b,l,\xi ^y(t,\omega )}
(u,\xi (u,\omega ))\circ $$
$$(I^{\otimes b}\otimes (a^y_{\xi ^y(t,\omega )})^{\otimes m}\otimes
(E^y_{\xi ^y(t,\omega )})^{\otimes l})|_{u=t} \} $$
due to Theorem $II.3.5$ \cite{lustpr} and Theorem $3.12$ above.
Therefore, $S_y(t,\tau ;\omega )\circ x\in U^1_y$ for each $x\in U^2_y$.
Glueing together local solutions with the help of transition functions
$\phi _{l,j}$ of charts with nonvoid intersections $U_l\cap U_j$
leads to the conclusion that a stochastic process $\xi $ is a solution
of the stochastic antiderivational equation
$(i)$ on $M$ if and only if for each
$t\in T$ for $\nu _{\xi (t)}$-a.e. $x\in M$ it coincides $P_{t,x}$-a.e.
with some local solution of this equation inside $U^2_x$, since
$\{ U^2_x: x\in M \} $ is a covering of $M$.
\par Consider a local representation $\xi ^{\phi }:=\phi (\xi )$,
then there exists the corresponding $S^{\phi }$ generated by
$d\xi ^{\phi }$ such that $S^{\phi }(t,\tau ;\omega )\circ \xi ^{\phi }
(\tau ,\omega )=\phi (S(t,\tau ;\omega )\circ \xi (\tau ,\omega ))$
for each $t$ and $\tau \in T\subset \bf K$, $\phi \in \{ \phi _j: j \} $.
\par In view of Proposition 3.23 there exists $\delta >0$ such that
$P \{ \omega : S_x(t,\tau ;\omega )\circ x \notin U^2_x \} \le P \{
\sup \| \phi (S_x(t,\tau ;\omega )\circ x) \| >1 \} \le C |t-\tau |$
for each $t, \tau \in T$ such that $|t-\tau |<\delta $,
where $C>0$ is a constant. Consider a family $\Upsilon $ of all
finite partitions $q$ of $T$ into disjoint unions of balls
$B({\bf K},t^q_k,r^q_k)$, where $t^q_k\in T$, $0<r^q_k\le \epsilon _q$,
$0<\epsilon _q<\delta $ for each $q\in \Upsilon $.
Let $q\le v$ if and only if $q\subset v$, then $\Upsilon $
is ordered by this relation. Consider a linearly ordered
subsequence $\Upsilon _0:=\{ q_k: k\in {\bf N} \} $
in $\Upsilon $ with $\lim_{k\to \infty }\sup \{ r^{q_k}_j: j \in q_k \} =0$
and for it define $\xi _k(t,\omega ):=S_{\xi _{k-1}
(t_l,\omega )} (t,t_l;\omega )\circ \xi _{k-1}(t_l,\omega )$
for each $k$ and $t\in B({\bf K},t^v_l,r^v_l)$ for each
$t_l\in q_{k-1}$ and each 
$k\ge 1$, where $v=q_{k-1}$, $\xi (t_0,\omega )=x$, $\xi _1
(t_l,\omega ):=\xi (t_l,\omega )$ for each $t_l\in q_1$.
Also define $S^k (t,t_0;\omega )\circ x=S_{\xi _{k-1}(t_l,\omega )}
(t,t_l;\omega )\circ \xi _{k-1}(t_l)$.
Consider $z(s,\omega ):=S_y(s,t^q_k;\omega )\circ y\in U^2_y$
for each $s\in B({\bf K},t^q_k,r^q_k)$.
For each $t\in B({\bf K},t^q_k,r^q_k)$ there is satisfied
the equality $S_y(t,t^q_k;\omega )\circ y=S_{z(s,\omega )}(t,s;\omega )
\circ z(s,\omega )$, since $S_y(t,t^q_k;\omega )\circ y=
S_y(t,s;\omega )\circ S_y(s,t^q_k;\omega )\circ y$ due to the existence
of a local solution.
\par Put $\Omega ^{\Upsilon _0}:=\bigcup_{k\in \Upsilon _0}\Omega _k$,
where $\Omega _k:=\bigcap_{l\in q_k}\Omega _{k,l}$, where
$\Omega _{k,l}:= \{ \omega : S_{\xi _k(t_l,\omega )}(s,t_l)\circ
\xi _k(t_l,\omega )\in U^2_{\xi _k(t_l,\omega )},
s\in B({\bf K},t^v_l,r^v_l) \} $.
From the existence of a local solution it follows, that
$S^k(t,t_0)\circ x=S^l(t,t_0)\circ x$ for each $k\ge l$ and each
$\omega \in \Omega _l$. In view of Theorem 3.19
there exists $\lim_{q_k\in \Upsilon _0}S^k(t,t_0;\omega )=S(t,t_0;\omega ).$
For each two linearly ordered subsets $\Lambda _1$ and $\Lambda _2$
in $\Upsilon $ there exists a linearly ordered subset $\Lambda $
in $\Upsilon $ such that $\Lambda \supset \Lambda _1\cup \Lambda _2$,
hence this limit does not depend on the choice of $\Upsilon _0$.
Events $\Omega _{k,l}$ and $\Omega _{k,j}$ are independent in total
for each $l\ne j$:
$P(\Omega _{k,l}\cap \Omega _{k,j})=P(\Omega _{k,l})P(\Omega _{k,j})$.
Since $\bf K$ is a finite algebraic extension
of the corresponding $\bf Q_p$, then there exists $n\in \bf N$ such that
$\bf K$ as the $\bf Q_p$-linear space is isomorphic with $\bf Q_p^n$.
Choose $\Upsilon _0$ such that for each $q_k\in \Upsilon _0$ the supremum
$\sup_{l\in q_k}r^{q_k}_l=:\delta _k\le p^{-k}$ and
$card (t_l: t_l\in q_k\cap B({\bf K},t_0,p^s))=:m_{k,s}\le p^{snk}$. 
In view of the ultrametric inequality from
$\| \alpha (\omega )+\beta (\omega ) \| \ge \delta $ it follows, that
$\max (|\alpha (\omega )|, |\beta (\omega )|)\ge \delta $
for each two random variables $\alpha $ and $\beta $. Therefore,
from Proposition 3.23 applied to $\phi _j(\xi ) -\phi _j(\xi _0)$,
and the inclusion $\xi (t,\omega )\in L^q(\Omega ,{\sf F},\lambda ;C^0(T,M))$
it follows, that $P \{ \Omega _k: t\in T\cap B({\bf K},t_0,p^s) \}
\ge (1 -C_kp^{-k})^{p^{snk}}$, where $\lim_kC_k=0$,
since $\xi (t,\omega )$ is uniformly continuous on $T\cap B({\bf K},t_0,p^s)$
for $\lambda $-a.e. $\omega $.
Therefore, $P \{ \Omega ^{\Upsilon _0}: t\in T\cap B({\bf K},t_0,p^s) \}
\ge \lim_k\exp (-C_ksn)=1$ for each given $s\in \bf N$.
\par From $S^k(t,t_0; \omega )\circ x=S^k(t,s;\omega )\circ
S^k(s,t_0;\omega )\circ x$ and taking the limit by $q\in \Upsilon $
it follows, that $S$ satisfies the evolution property
$S(t,t_0;\omega )\circ x=S(t,s;\omega )\circ S(s,t_0;\omega )\circ x$.
Then $S(t,t_0;\omega )\circ x$ is a measurable function
of $x$, since it is a superposition $S(t,t_0;\omega )\circ x=
S^k(t,t_0;\omega )\circ x$ of locally measurable functions.
\par {\bf 3.25. Corollary.} {\it Let conditions of Theorem 3.24
be satisfied, $f$ be a function on $T\times M$ such that
each $f\circ \phi _j^{-1}$ satisfies Conditions $3.11.(iii,iv)$
on its domain. Then a generating operator of an evolution 
family $S(t,\tau ;\omega )$ of a stochastic process 
$\eta (t,\omega )=f(t,\xi (t,\omega ))$ is given by the equation:
$$(i)\quad A(t;\omega )\eta (t,\omega )=
\sum_{m+b\ge 1, 0\le m\in {\bf Z}, 0\le b\in {\bf Z}}
((m+b)!)^{-1}\sum_{l=0}^m{{m+b}\choose m}$$
$$\sum_{l_1+...+l_m=m-l,n_1+...+n_m=l} \{
b({\hat P}_{u^{b+m-l-1},w(u,\omega )^l}[(\partial ^b\nabla ^m f/
\partial u^b\partial x^m) (u,\xi (u,\omega ))\circ
(a_{l_1,n_1}\otimes $$
$$... \otimes a_{l_m,n_m})\circ (I^{\otimes (b-1)}\otimes
a^{\otimes (m-l)}\otimes E^{\otimes l})])|_{u=t}+$$
$$(m-l)({\hat P}_{u^{b+m-l-1},w(u,\omega )^l}[(\partial ^b \nabla ^mf/
\partial u^b\partial x^m) (u,\xi (u,\omega ))\circ
(a_{l_1,n_1}\otimes ... \otimes a_{l_m,n_m})\circ $$
$$(I^{\otimes b}\otimes
a^{\otimes (m-l-1)}\otimes E^{\otimes l})]a)|_{u=t}+$$
$$ l({\hat P}_{u^{b+m-l},w(u,\omega )^{l-1}}[(\partial ^b\nabla ^mf/
\partial u^b\partial x^m) (u,\xi (u,\omega ))\circ
(a_{l_1,n_1}\otimes ... \otimes a_{l_m,n_m})\circ $$
$$(I^{\otimes b}\otimes
a^{\otimes (m-l)}\otimes E^{\otimes (l-1)})]E{w'}_u(u,\omega ))
|_{u=t} \} .$$}
\par {\bf Proof.} In view of Theorem $3.24$
there exists a generating operator $S(t,\tau ;\omega )$
of an evolution family. For each chart $(U_j,\phi _j)$
the stochastic process $f\circ \phi _j^{-1}(\xi )$ is given by
Equation $3.12.(i)$. Consider the covariant differentiation
$(\nabla f/\partial x).h=\nabla _hf$ on the manifold $M$,
where $h\in T_xM$. For a random variable belonging to
$L^q(\Omega ,{\sf F},\lambda ;C^1(M,X))$ its derivative
and partial difference quotients ${\bar \Phi }^1f\circ \phi _j^{-1}(x;h;b)$
are naturally understood as elements of the corresponding
spaces $L^q(\Omega ,{\sf F},\lambda ;C^0(W_j,X))$ such that
each limit $\lim_{x\to x_0}g(x,\omega )=c(\omega )$ is taken
in $L^q(\Omega ,{\sf F},\lambda ;C^0(M,X))$, where $W_j:=
\{ (x,h,b)\in U_j\times X\times {\bf K}: x+bh\in U_j \} $.
In another words it exists if and only if
$\lim_{x\to x_0} \| g(x,\omega )-c(\omega ) \| _{L^q}=0$,
where $c\in L^q(\Omega ,{\sf F},\lambda ;X)$.
Then $f(t,\xi (t,\omega ))=\lim_k f(t,S^k (t,t_0;\omega )\circ x).$
For each chart put $f_j(t,*):=f(t,\phi _j^{-1}(*)\phi _j\circ S_y)$,
where $S_y(t,\tau ;\omega )y$ does not leave for $\lambda $-a.e.
$\omega \in \Omega $ a clopen subset $U_j$ in $M$
for each $t$ and $\tau \in T_j$, $T_j\subset T$, $\bigcup_jT_j=T$,
$T_j$ is a clopen subset in $T$. Then define
$(f(t,S^k(t,\tau ;\omega )x)'_x.h=f'_x(t,S^k(t,\tau )\circ x)
A^k(t,\tau ;\omega )h$ and take the limit.
From Lemma $I.2.3$ \cite{lustpr} it follows the statement of this Corollary.
\par {\bf 3.26. Remarks.} In the particular case $3.10$
Formula $3.25.(i)$ simplifies. When the family of $\Gamma $ together
with all its covariant derivatives along $a$ and $Ew'$
is equiuniformly bounded on each $U_j$, then
$3.25.(i)$ can be written in another form using the identity
$\nabla ^{m+1}f.(h_1,...,h_m)=\nabla _{h_{m+1}}(\nabla ^m
f.(h_1,...,h_m))-\sum_{l=1}^m(\nabla ^mf).(h_1,...,h_{l-1},
\nabla _{h_{m+1}}h_l,h_{l+1},...,h_m)$, in particular with
$h_l\in \{ a, Ew' \} $.
\par In a weak topology $C^0(T,c_0(\alpha ,{\bf K}))$
is isomorphic with $c_0(\alpha ,{\bf K})^T$. Let $\theta : {\bf K}\to
\bf R$ be a continuous surjective quotient mapping such that
$\theta (B({\bf K},0,1))=[0,1]$. Then for each $\xi \in
L^q(\Omega ,{\sf F},\lambda ;c_0(\alpha ,{\bf K})^T)$ there exists
$\theta (\xi )\in L^q(\Omega ,{\sf F},\lambda ;$ \\
$c_0(\alpha ,{\bf R})^{\theta (T)})$, that induces a surjective mapping
$\theta ^*$ from $L^q(\Omega ,{\sf F},\lambda ;c_0(\alpha ,{\bf K})^T)$
onto $L^q(\Omega ,{\sf F},\lambda ;c_0(\alpha ,{\bf R})^{\theta (T)})$.
Therefore, for each stochastic process $\eta $ in \\
$L^q(\Omega ,{\sf F},\lambda ;c_0(\alpha ,{\bf R})^{\theta (T)})$
there exists a stochastic process $\xi $ in
$L^q(\Omega ,{\sf F},\lambda ;c_0(\alpha ,{\bf K})^T)$
such that $\theta (\xi )=\eta $. On the other hand,
$\bf K$ is a projective limit of discrete rings $S_{\pi ^n}$
isomorphic to the quotient of $\bf K$ by the equivalence relation
associated with disjoint subsets $x_j+\pi ^nB({\bf K},0,1)$ in $\bf K$,
$j=0,1,2,...$, $x_0:=0$, $\pi \in \bf K$, $|\pi |=\max \{ |x|:
|x|<1, x\in {\bf K} \} $. Therefore,
$L^q(\Omega ,{\sf F},\lambda ;c_0(\alpha ,{\bf K})^T)$
is isomorphic as the topological space to projective limit of modules
$L^q(\Omega ,{\sf F},\lambda ;c_0(\alpha ,{\bf S_{\pi ^n}}
)^{\bf S_{\pi ^n}})$ over discrete rings $\bf S_{\pi ^n}$,
since simple functions
are dense in $L^q$, consequently, $\xi $ is equal to the projective
limit of stochastic processes with values in discrete modules
$c_0(\alpha ,{\bf S_{\pi ^n}})$ over rings $\bf S_{\pi ^n}$.
This opens a possibility of approximation of stochastic processes
by stochastic processes with values in discrete modules.
Certainly there is not any simple relation between classical
and non-Archimedean stochastic equations, so if $\xi $ satisfies
definite stochastic antiderivational equation relative to $w$
it is a problem to find a classical stochastic equation to which
$\theta (\xi )$ satisfies relative to either $\theta (w)$ or
a standard stochastic process (Wiener, L\`evy) and vice versa.
\par Theorem 3.24 and Corollary 3.25 are applicable in particular to
totally disconnected Lie groups over non-Archimedean
fields of characteristic zero.

{\underline {Acknowledgment}}. The author thanks
Prof. J. Oesterle and Prof. H. Rosenberg for hospitality
at the "Institut de Math\'ematiques Universit\'e Paris 6 et 7" and
the "Ministere de la Recherche" and "\'EGIDE" of France for support.
\par Address: Theoretical Department, Institute of General Physics,
\par Russian Academy of Sciences,
\par Str. Vavilov 38, Moscow, 119991 GSP-1, Russia
\par E-mail: ludkovsk@fpl.gpi.ru

\begin{thebibliography}{99}
\bibitem{beldal} Ya. I. Belopolskaya, Yu. L. Dalecky.
"Stochastic equations and differential geometry"
(Dordrecht: Kluwer, 1989).
\bibitem{bikvol} A.H. Bikulov, I.V. Volovich.
"$p$-adic Brownian motion". Izv. Ross. Akad. Nauk.
Ser. Math. {\bf 61: 3} (1997), 75-90.
\bibitem{boui}
N. Bourbaki. "Integration". Chapters 1-9 (Moscow: Nauka, 1970 and 1977).
\bibitem{coxmil} D.R. Cox, H.D. Miller. "The theory of stochastic processes"
(London: Chapman and Hall, 1995).
\bibitem{dalfo} Yu. L. Dalecky, S.V. Fomin. "Measures and differential 
equations in infinite dimensional space" (Dordrecht: Kluwer, 1991).
\bibitem{eng} R. Engelking. "General topology" (Moscow: Mir, 1986).
\bibitem{evans}  S.N. Evans. "Continuity properties
of Gaussian stochastic processes indexed by a local field".
Proceed. Lond. Math. Soc. Ser. 3, {\bf 56} (1988), 380-416.
\bibitem{freput} J. Fresnel, M. van der Put. "G\'eom\'etrie
analytique rigide et applications" (Boston, Birkh\"auser, 1981).
\bibitem{gihsko} I.I. Gihman, A.V. Skorohod. "Stochastic differential
equations and their apllications" (Kiev: Naukova Dumka, 1982).
\bibitem{gisk3} I.I. Gihman, A.V. Skorohod. "Theory of stochastic processes".
{\bf V.V. 1-3} (Moscow: Nauka, 1975).
\bibitem{hew} E. Hewitt, A. Ross. "Abstract harmonic analysis"
(Berlin: Springer, 1979).
\bibitem{ikwat} N. Ikeda, S. Watanabe. "Stochastic differential equations
and diffusion processes" (Moscow: Nauka, 1986).
\bibitem{itkean} K. It$\hat o$, H.P. McKean. "Diffusion
processes and their sample paths" (Berlin: Springer, 1996).
\bibitem{jang} Y. Jang. "Non-Archimedean quantum mechanics".
Tohoku Mathem. Publications. $N^o$ {\bf 10} 
(Tohoku: Toh. Univ., Math. Inst., 1998).
\bibitem{khrif} A.Yu. Khrennikov. "Generalized functions and Gaussian 
path integrals". Izv. Acad. Nauk. Ser. Mat. {\bf 55} 
(1991), 780-814.
\bibitem{khipr} A. Khrennikov. "Interpretations of probability"
(Utrecht: VSP, 1999).
\bibitem{kling} W. Klingenberg. "Riemannian geometry"
(Berlin: Walter de Gruyter, 1982).
\bibitem{kolfp} A.N. Kolmogorov. "Foundations of the theory of probability"
(New York: Chelsea Pub. Comp., 1956).
\bibitem{luumnls} S.V. Ludkovsky. "Quasi-invariant measures on
non-Archimedean semigroups of loops". Russ. Math. Surveys {\bf 53: 3} (1998),
633-634.
\bibitem{luseamb} S.V. Ludkovsky. "Irreducible unitary representations
of non-Archimedean groups of diffeomorphisms". Southeast Asian
Mathem. Bull. {\bf 22: 3} (1998), 301-319.
\bibitem{luseamb2} S.V. Ludkovsky.
"Poisson measures for topological groups and their representations".
Southeast Asian Mathem. Bull. {\bf 25: 4} (2002), 653-680.
\bibitem{lufpm} S.V. Ludkovsky.
"Non-Archimedean polyhedral decompositions of
ultrauniform spaces". Fundam. i Prikl. Math. {\bf 6: 2} (2000), 455-475. 
\bibitem{lubp99} S.V. Ludkovsky. "Properties of quasi-invariant measures on
topological groups and associated algebras".
Annales Math\'ematiques Blaise Pascal. {\bf 6: 1} (1999), 33-45.
\bibitem{lutmf99} S.V. Ludkovsky. "Measures on groups of diffeomorphisms
of non-Archimedean manifolds, representations of groups and their 
applications". Theoret. and Math. Phys. {\bf 119: 3} (1999), 698-711.
\bibitem{lubp2} S.V. Ludkovsky. "Quasi-invariant measures on non-Archimedean 
groups and semigroups of loops and paths, their representations. I, II".
Annales Math\'ematiques Blaise Pascal. {\bf 7: 2} (2000), 19-53, 55-80.
\bibitem{lulapm} S.V. Ludkovsky. "Quasi-invariant and 
pseudo-differentiable measures
on a non-Archimedean Banach space. I, II". Los Alamos Preprints
{\bf math.GM/0106169} and {\bf math.GM/0106170} (http://xxx.lanl.gov/;
earlier version: ICTP {\bf IC/96/210},  October 1996, 50 p.p.;
http://www.ictp.trieste.it/).
In part it will appear in Analys. Mathem. and in J. Mathem. Sci.
\bibitem{lufpmsp} S.V. Ludkovsky. "Stochastic processes
on groups of diffeomorphisms and loops of real, complex and
non-Archimedean manifolds". Fundam. i Prikl. Mathem.
{\bf 7: 4} (2001), 1091-1105.
\bibitem{lustpr} S.V. Ludkovsky.
"Stochastic processes on non-Archimedean spaces". "I.
Stochastic processes on Banach spaces".
"II. Stochastic antiderivational equations".
"III. Stochastic processes on totally disconnected topological groups".
Intern. J. of Math. and Math. Sci., is accepted to publication
(previous variant: Los Alamos National Laboratory, USA.
Preprints {\bf math.GM/0104069}, {\bf math.GM/0104070},
{\bf math.GM/0106132}, April-June 2001).
\bibitem{lukhr} S.V. Ludkovsky, A. Khrennikov. "Stochastic processes
on non-Archimedean spaces with values in non-Archimedean
fields". Markov Processes and Related Fields. {\bf 8} (2002), 1-34
(earlier version: Los Alamos National Laboratory, USA.
Preprint {\bf math.CA/0110305}, 20 pages, 06 April 2001).
\bibitem{malb} P. Malliavin. "Stochastic analysis" (Berlin: Springer, 1997).
\bibitem{mckean} H.P. Mc Kean. "Stochastic integrals" 
(Moscow: Mir, 1972).
\bibitem{nari} L. Narici, E. Beckenstein. "Topological vector spaces"
(New York: Marcel-Dekker Inc., 1985).
\bibitem{oksb} B. $\O$ksendal. "Stochastic differential equations"
(Berlin: Springer, 1995).
\bibitem{pietsch} A. Pietsch. "Nukleare lokalkonvexe R\"aume"
(Berlin: Akademie-Verlag, 1965).
\bibitem{roo} A.C.M. van Rooij. "Non-Archimedean fucntional analysis"
(New York: Marcel Dekker, 1978). Ser. Pure and Appl. Math. V. {\bf 51}.
\bibitem{sch1} W.H. Schikhof. "Ultrametric calculus" (Cambridge:
Cambr. Univ. Press, 1984).  
\bibitem{vla3} V.S. Vladimirov, I.V. Volovich, E.I. Zelenov. "$p$-Adic
analysis and mathematical physics" (Moscow: Fiz.-Mat. Lit, 1994).
\bibitem{wei} A. Weil. "Basic number theory" (Berlin: Springer, 1973).
\end{thebibliography}
\end{document}